# Optimal *a priori* balance in the design of controlled experiments

Nathan Kallus

*Cornell University, New York, USA*

**Summary.** We develop a unified theory of designs for controlled experiments that balance baseline covariates *a priori* (before treatment and before randomization) using the framework of minimax variance and a new method called kernel allocation. We show that any notion of *a priori* balance must go hand in hand with a notion of structure, since with no structure on the dependence of outcomes on baseline covariates complete randomization (no special covariate balance) is always minimax optimal. Restricting the structure of dependence, either parametrically or non-parametrically, gives rise to certain covariate imbalance metrics and optimal designs. This recovers many popular imbalance metrics and designs previously developed *ad hoc*, including randomized block designs, pairwise-matched allocation and rerandomization. We develop a new design method called kernel allocation based on the optimal design when structure is expressed by using kernels, which can be parametric or non-parametric. Relying on modern optimization methods, kernel allocation, which ensures nearly perfect covariate balance without biasing estimates under model misspecification, offers sizable advantages in precision and power as demonstrated in a range of real and synthetic examples. We provide strong theoretical guarantees on variance, consistency and rates of convergence and develop special algorithms for design and hypothesis testing.

## 1. Introduction

Achieving balance between experimental groups is a cornerstone of causal inference; otherwise any observed difference may be attributed to a difference other than the treatment alone. In clinical trials, and more generally controlled experiments, the experimenter controls the administration of treatment and complete randomization of subjects has been the gold standard for achieving this balance.

The expediency and even necessity of complete randomization, however, have been controversial since the founding of statistical inference in controlled experiments. William Gosset, 'Student' of Student's *t*-test, said of assigning field plots to agricultural interventions that it 'would be pedantic to continue with an arrangement of [field] plots known beforehand to be likely to lead to a misleading conclusion', such as arrangements in which one experimental group is on average higher on a 'fertility slope' than the other experimental group (Student, 1938). Of course, as the opposite is just as likely under complete randomization, this is not an issue of estimation bias in its modern definition, but of estimation variance. Gosset's sentiment is echoed in the common statistical maxim 'block what you can; randomize what you cannot' that has been attributed to George Box and in the words of such individuals as Heckman (2008) ('Randomization is a metaphor and not an ideal or "gold standard"') and Rubin (2008) ('For gold standard answers, complete randomization may not be good enough'). In one interpretation, these can be seen as calls for the experimenter to ensure that experimental groups are balanced at the outset of the experiment, before applying treatments and before randomization.

There is a variety of designs for controlled experiments that attempt to achieve better balance in terms of measurements made before treatment, known as baseline covariates, under the understanding that a predictive relationship possibly holds between baseline covariates and the outcomes of treatment. We term this sort of approach *a priori balancing* as it is done before applying treatments and before randomization (the term *a priori* is chosen to contrast with *post hoc* methods such as poststratification, which may be applied after randomization and after treatment as in McHugh and Matts (1983)). The most notable *a priori* balancing designs are randomized block designs (Fisher, 1935), pairwise-matched allocation (Greevy *et al.*, 2004) and rerandomization (Morgan and Rubin, 2012). There are also sequential methods for allocation that must be decided as subjects are admitted but before treatment is applied (Efron, 1971; Pocock and Simon, 1975; Kapelner and Krieger, 2014).

Each of these implicitly defines imbalance between experimental groups differently. Blocking attempts to achieve exact matching (when possible): a binary measure of imbalance that is 0 only if the experimental groups are identical in their discrete or coarsened baseline covariates. Pairwise-matched allocation treats imbalance as the sum of pairwise distances, given some pairwise distance metric such as the Mahalanobis distance. There are both globally optimal and greedy heuristic methods that address this imbalance measure (Gu and Rosenbaum, 1993). Morgan and Rubin (2012) defined imbalance as the Mahalanobis distance in group means and proposed rerandomization as a heuristic method for reducing it (non-optimally). Bertsimas *et al.* (2015) defined imbalance as the sum of discrepancies in group means and variances and found optimal designs with globally minimal imbalance.

It is not immediately clear when each of these different characterizations of imbalance is appropriate and when, for that matter, deviating from complete randomization is justified. Moreover, the connection between imbalance that is measured before treatment and the estimation error after treatment is often unclear. We here argue that, without structural information on the dependence of outcomes on baseline covariates, there need not be any such connection and complete randomization is minimax optimal—a 'no-free-lunch' theorem (Wolpert and Macready, 1997). However, when structural knowledge is expressed as membership of conditional expectations in a normed vector space of functions (specifically, Banach spaces), different minimax optimal rules arise for the *a priori* balancing of experimental groups, based on controlling after-treatment variance via a before-treatment imbalance metric. We show how certain choices of such structure reconstruct each of the aforementioned methods or associated imbalance metrics. We also develop and study a new range of powerful methods, kernel allocation, that arises from choices of structure, both parametric and non-parametric, based on kernels. Theoretical and numerical evidence suggests that kernel allocation should be favoured in practice.

We study in generality the characteristics of any method that arises from our framework, including its estimation variance and consistency, intimately connecting *a priori* balance to post-treatment estimation. Whenever a parametric model of dependence is known to hold, we show that, relative to complete randomization, the variance due to the optimal design converges linearly ($2^{-\Omega(n)}$ for $n$ subjects) to the best theoretically possible variance, achieving nearly perfect covariate balance without biasing estimates in case of misspecification. This both generalizes and formalizes an observation that was made by Bertsimas *et al.* (2015). We develop algorithms for computing optimal designs by using mixed integer programming and semidefinite programming and provide hypothesis tests for these designs. A Python package implementing all of these algorithms is provided at `http://www.nathankallus.com` and

> `http://wileyonlinelibrary.com/journal/rss-datasets`

## 2.   The effect of structural information and lack thereof

We begin by describing the set-up. Let $m$ denote the number of treatments to be investigated (including controls). We index the subjects by $i = 1, \ldots, n$ and assume that $n = mp$ is divisible by $m$. We assume that the subjects are independently randomly sampled but we shall consider estimating both sample and population effects. We denote assigning subject $i$ to a treatment $k$ by $W_i = k$. We let $w_{ik} = \mathbb{I}(W_i = k)$ and $W = (W_1, \ldots, W_n)$. When $m = 2$, we shall use $u_i = w_{i1} - w_{i2}$. The design is the distribution over $W$ induced by the researcher. An assignment is a single realization of $W$.

As is common for controlled trials, we assume non-interference (see Cox (1958), Rubin (1986) and Rosenbaum (2007)), i.e. a subject who is assigned to a certain treatment exhibits the same outcome regardless of others' assignments. Under this assumption we can define the (real-valued) potential post-treatment outcome $Y_{ik} \in \mathbb{R}$ of subject $i$ if it were to be subjected to the treatment $k$. We let $Y = ((Y_{11}, \ldots, Y_{1m}), \ldots, (Y_{n1}, \ldots, Y_{nm}))$ denote the collection of all potential outcomes of all $n$ subjects. We assume throughout that each $Y_{ik}$ has second moments. Let $X_i$, taking values in some $\mathcal{X}$, be the baseline covariates of subject $i$ that are recorded before treatment and let $X = (X_1, \ldots, X_n)$. The space $\mathcal{X}$ is general; assumptions about it will be specified as necessary. As an example, it can be composed of real-valued vectors $\mathcal{X} \subset \mathbb{R}^d$ that include both discrete (dummy) and continuous variables.

We denote by $\text{TE}_{kk'i} = Y_{ik} - Y_{ik'}$ the unobservable causal treatment effect for subject $i$. Two unobservable quantities of estimation interest are the *sample average (causal) treatment effect*, SATE,

$$\text{SATE}_{kk'} = \frac{1}{n} \sum_{i=1}^{n} \text{TE}_{kk'i} = \frac{1}{n} \sum_{i=1}^{n} Y_{ik} - \frac{1}{n} \sum_{i=1}^{n} Y_{ik'},$$

and the *population average (causal) treatment effect*, PATE,

$$\text{PATE}_{kk'} = \mathbb{E}[\text{TE}_{kk'1}] = \mathbb{E}[\text{SATE}_{kk'}].$$

When subjects are randomly sampled from the population, SATE is an unbiased and consistent estimate of PATE. Our estimator for either SATE or PATE will always be the simple mean differences estimator

$$\hat{\tau}_{kk'} = \frac{\sum\limits_{i : W_i = k} Y_{ik}}{\sum\limits_{i : W_i = k} 1} - \frac{\sum\limits_{i : W_i = k'} Y_{ik'}}{\sum\limits_{i : W_i = k'} 1}.$$

Note that, as an estimator for PATE, the estimator $\hat{\tau}_{kk'}$ precludes *post hoc* adjustments by regression (Freedman, 2008). Without parametric assumptions that are necessary for such adjustments and given a random design $\hat{\tau}_{kk'}$ is the uniformly minimum variance unbiased estimator for PATE (Lehmann and Casella (1998), chapter 2). Although carefully conducted *post hoc* regression adjustments need only to introduce slight bias in misspecified settings (Lin, 2013), as we shall see, optimal *a priori* balance can achieve the same efficiency with zero finite sample bias. We drop subscripts when $m = 2$ and set $k = 1$ and $k' = 2$.

Throughout we shall consider only designs that satisfy the following conditions.

*Assumption 1.* The design

(a)  does not depend on future information, i.e. $W$ is independent of $Y$, conditionally on $X$,
(b)  blinds (randomizes) the identity of treatments, i.e., for any permutation $\pi$ of $1, \ldots, m$, we have that $\mathbb{P}\{W = (k_1, \ldots, k_n) | X\} = \mathbb{P}\{W = (\pi(k_1), \ldots, \pi(k_n)) | X\}$ and
(c)  splits the sample evenly, i.e. surely $\Sigma_{i:W_i=k} 1 = p \ \forall k$.

We interpret assumption 1 as the definition of *a priori balance* as it requires that all balancing be done before applying treatments (condition (a)), and before randomization (conditions (b) and (c)). Condition (a) is a reflection of the temporal logic of first assigning, and then experimenting.

Condition (b) says that balancing is done before randomization. This ensures that the estimators $\hat{\tau}_{kk'}$ resulting from the design are always unbiased, both conditionally on $X$ and $Y$ (i.e. in estimating SATE) and marginally (i.e. in estimating PATE; more detail is given in theorem 7 in Section 3). Condition (b) corresponds to the blinding aspect of randomization, originally described by Fisher (see Senn (1994)).

Condition (c) is one way to achieve condition (b) in non-completely randomized designs. If $W$ is an even assignment then randomly permuting treatment indices will blind their identity. Otherwise, fixing an uneven assignment, a treatment can be identified by the size of its experimental group.

Assumption 1 defines the space of designs that we shall be comparing. When optimizing balance, it will be over these designs. Further, we consider randomized designs. Given fixed subjects, we denote by $\sigma(W)$ the probability of assignment $W$ where $\sigma$ is the distribution over assignments that are prescribed by the design. We denote by $\mathcal{W} \subset \{1, \ldots, m\}^n$ the space of feasible assignments satisfying condition (c) and by $\Delta \subset [0, 1]^{|\mathcal{W}|}$ the space of feasible distributions over assignments satisfying all of assumption 1, i.e. all *a priori* balancing designs given the fixed subjects.

## 2.1. 'No free lunch'

We shall now argue that, without structural information on the relationship between $X_i$ and $Y_{ik}$, complete randomization is minimax optimal. A minimax characterization of optimality assumes an adversarial nature because, in the absence of information, uncertainty can neither be inscribed nor distributionally described. One of Fisher's justifications for randomization is as a blinding instrument to achieve unbiasedness in the face of uncertainty (Senn, 1994). Beyond blinding and unbiasedness (i.e. assumption 1, part (b)), randomization also reduces variance in the face of uncertainty. The type of variance that we study here is that due to the design, while experimental units are fixed but outcomes unknown (i.e. variance conditioned on $X$ and $Y$). A different type of variance is that due to a Bayesian distributional characterization of unknown outcomes (not fixed) where any variance due to the design is redundant because in a basic application of Bayesian optimality a single assignment is always optimal, losing the guarantee of unbiasedness in the face of uncertain outcomes. A full discussion of the relevance of a worst-case analysis and a comparison with a Bayesian approach is given in Section 2.4.1.

The results herein show that, for any design satisfying the general conditions of assumption 1, there is some value for outcomes such that the variance due to the randomization of the design is never better than complete randomization. We restrict here to $m = 2$.

Among estimators that are unbiased, the standard way of comparing efficiency is variance. By the law of total variance and by conditional unbiasedness under a design satisfying assumption 1, we have

$$\mathrm{var}(\hat{\tau}) = \mathbb{E}[\mathrm{var}(\hat{\tau}|X, Y)] + \mathrm{var}(\mathrm{SATE}).$$

The variance of SATE is independent of our choice of *a priori* balancing design since SATE is equal to the difference of all outcomes of all subjects regardless of assignment—in other words, the definition of SATE includes only $Y$ (not $W$) and assumption 1, part (a), ensures independence. The choice of design can affect only the first term. Therefore, the efficient designs that satisfy assumption 1 are exactly those that minimize the conditional variance $\mathrm{var}(\hat{\tau}|X,Y)$ for the given realized subjects, i.e. those that optimally treat each realization of subjects separately, on a sample path by sample path basis.

However, in practice, we do not know $Y$: only $X$ (assumption 1, part (a)). Since we assume no structural information on their relationship, we consider an adversarial nature that chooses $Y$ so as to increase our variance. To avoid a trivial infinity, we fix the variance of the complete-randomization estimator (denoted $\hat{\tau}^{\mathrm{CR}}$), which is the same as fixing the overall magnitude of $Y$. The following theorem shows that, in this situation, complete randomization is optimal.

*Theorem 1.* Fix $X \in \mathcal{X}^n$, $V > 0$ and $m = 2$. Then, among designs satisfying assumption 1, complete randomization minimizes

$$\max_{Y \in \mathbb{R}^{n \times m} : \mathrm{var}(\hat{\tau}^{\mathrm{CR}}|X, Y=V)} \mathrm{var}(\hat{\tau}|X, Y),$$

More generally, for any row permutationally invariant seminorm $\|\cdot\|$ on $\mathbb{R}^{n \times m}$ (i.e. $\|PY\| = \|Y\|$ for all $n \times n$ permutation matrices $P$), complete randomization minimizes

$$\max_{Y \in \mathbb{R}^{n \times m} : \|Y\| \leqslant V} \mathrm{var}(\hat{\tau}|X, Y).$$

In particular, if a design fixes a partition of units and only flips a coin for blinding treatments,

$$\max_{Y \in \mathbb{R}^{n \times m} : \mathrm{var}(\hat{\tau}^{\mathrm{CR}}|X, Y)=V} \mathrm{var}(\hat{\tau}|X, Y) = V(n-1). \tag{2.1}$$

### 2.1.1. Example 1

To showcase the result of theorem 1, we now construct an example in the style of Cochran and Cox (1957), section 4.26, where we fix values for outcomes and study the variance that is due only to the randomization of various designs. Fix $n = 2^b$ a power of 2, $m = 2$ and any $\tau > 0$. Let

$$X_i = \sum_{t=0}^{b-\max\{2, \log_2(i)\}} \left\{ (-1)^{\lceil i/2^{t-1} \rceil} \times 2^{-2^{b-1} + 2^{b-t-1} + (i-1 \bmod 2^{t-1})} \right\},$$

$$Y_{i1} = (-1)^i = (-1)^{\log_2\{\mathrm{round}(|X_i|)\}} - \tau/2,$$

$$Y_{i2} = Y_{i1} + \tau.$$

This rather complicated construction essentially yields

$$X \approx \mathrm{round}(X) = (-1, -2, -4, \ldots, -2^{2^{b-1}-1}, 1, 2, 4, \ldots, 2^{2^{b-1}-1})$$

with just enough perturbation so that the assignment $W = (1, 2, 1, 2, \ldots, 1, 2)$ *uniquely* minimizes Mahalanobis distance between group means. For blocking, we block $X$ into eight consecutive intervals so that each contains the same number of subjects: $2^{b-3}$ (for $b \geqslant 4$). For pairwise-matched allocation, we use optimal matching with the pairwise Mahalanobis distance. And, for rerandomizaiton, we use the procedure of Morgan and Rubin (2012) with an infinitesimal acceptance probability that essentially globally minimizes the Mahalanobis distance between group means. We plot the resulting conditional variances $\mathrm{var}(\hat{\tau}|X, Y)$ in Fig. 1. Complete randomization has variance $4/(n-1)$, blocking has $4/(n-8)$ (in agreement with Cochran and

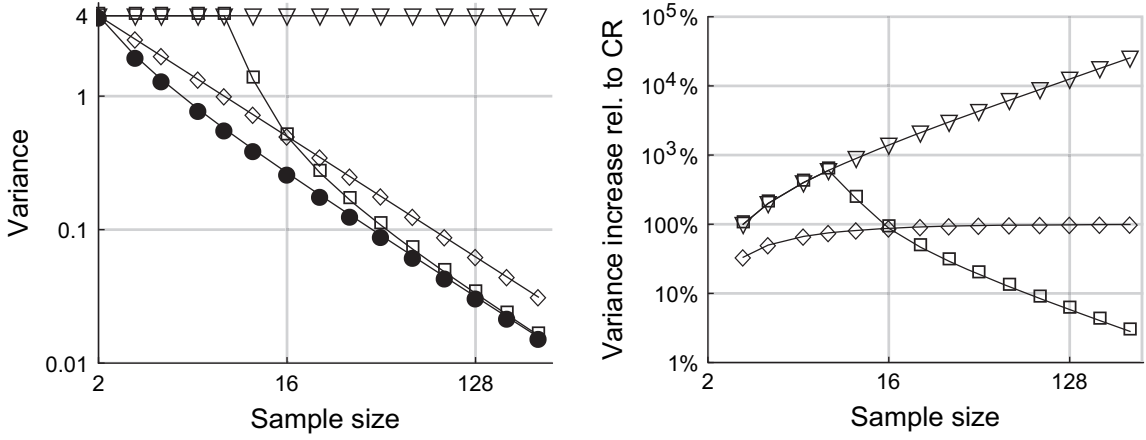

**Fig. 1.** Estimation variance under various designs in example 1 (note the logarithmic scales): ●, complete randomization; ▽, rerandomization $p \ll 1$; □, blocking; ◇, pairwise matching

Cox (1957), equation (4.3)), pairwise-matched allocation has $8/n$ and rerandomization with infinitesimal acceptance probability has 4, realizing equation (2.1) exactly.

In this example, complete randomization always does better than various standard designs, providing a concrete example of the conclusion of theorem 1. The construction, however, is plainly perceived as contrived. But, to characterize this construction as contrived, we must have a notion of structure that explains what is uncontrived. Therefore, a notion of balance must go hand in hand with a notion of structure. Next, we propose one way to formulate this.

### 2.2. Structural information and optimal designs

We now consider the effect of restricting nature by imposing particular structure on the conditional expectations of outcomes. Denote

$$f_k(x) := \mathbb{E}[Y_{ik}|X_i = x]$$

and

$$\epsilon_{ik} := Y_{ik} - f_k(X_i).$$

The non-random function $f_k$ is called the *conditional expectation function*. The law of iterated expectation yields that $\epsilon_{ik}$ has mean 0 and is mean independent of $X_i$. Combined with independence of subjects and equal sample size among treatments, this yields (see theorem 7 in Section 3)

$$\operatorname{var}(\hat{\tau}) = \mathbb{E}[\operatorname{var}\{B(W, \hat{f})|X\}] + \frac{1}{n}\operatorname{var}(\epsilon_{11} + \epsilon_{12}) + \operatorname{var}(\text{SATE}),$$

$$B(W, f) = \frac{2}{n}\sum_{i: W_i = 1} f(X_i) - \frac{2}{n}\sum_{i: W_i = 2} f(X_i), \quad \hat{f}(x) = \frac{f_1(x) + f_2(x)}{2}. \tag{2.2}$$

Throughout, we shall use $f : X \to \mathbb{R}$ to mean any generic function. By assumption 2, part (b), $\mathbb{E}[B(W, f)] = 0$ for any $f$ and hence $\operatorname{var}\{B(W, f)|X\} = \mathbb{E}[B(W, f)^2|X]$. As before, the variances of SATE and of $\epsilon_{11} + \epsilon_{12}$ are independent of our choice of design. Therefore, efficient designs satisfying assumption 1 are exactly those that minimize $\mathbb{E}[B(W, \hat{f})^2|X]$ for the given realized baseline covariates $X$, i.e. an efficient design optimally treats each realization of $X$ separately.

The unknown that determines which design is efficient is reduced to $\hat{f}$. We let nature choose it adversarially, fixing only the magnitude of $\hat{f}$. To define a magnitude of $\hat{f}$, we assume that $f_k \in \mathcal{F}$ $\forall\, k$, where $\mathcal{F}$ is a normed vector space with norm $\|\cdot\| : \mathcal{F} \to \mathbb{R}_+$. This will represent our structural information about the dependence between $X_i$ and $Y_{ik}$. This space is a subspace of the vector space $\mathcal{V}$ of all functions $\mathcal{X} \to \mathbb{R}$ under pointwise addition and scaling. For functions $f$ that are not in $\mathcal{F}$ we formally define $\|f\| = \infty$ so that we have a magnitude function $\|\cdot\| : \mathcal{V} \to \mathbb{R} \cup \{\infty\}$. Thus, when $\mathcal{F}$ is finite dimensional, the assumption that $\|f_k\| < \infty$ is the assumption of a parametric model.

In the context of minimax variance as the unknown $\hat{f}$ varies, with all structural information summarized by $\|\hat{f}\| < \infty$, we are interested in minimizing

$$\max_{\|f\| \leqslant \|\hat{f}\|} \mathbb{E}[B^2(W, f)|X] = \|\hat{f}\|^2 \max_{\|f\| \leqslant 1} \mathbb{E}[B^2(W, f)|X] = \|\hat{f}\|^2 \max_{f \in \mathcal{F}} \frac{\mathbb{E}[B^2(W, f)|X]}{\|f\|^2}, \qquad (2.3)$$

where the equalities hold because $B(W, \alpha f) = \alpha B(W, f)$ and $\|\alpha f\| = |\alpha| \|f\|$. Hence, which are the minimizers of the above equation is completely invariant to the specific value of $\|\hat{f}\| < \infty$ and all that matters is that it is in fact finite, or equivalently that $\hat{f} \in \mathcal{F}$.

Note that $B(W, f)$ is invariant to constant shifts to $f$, i.e. $B(W, f) = B(W, f + c)$ where $c \in \mathbb{R}$ represents a constant function $x \mapsto c$. To factor this artefact away, we consider the quotient space $\mathcal{F}/\mathbb{R}$, which consists of elements $[f] = \{f + c : c \in \mathbb{R}\}$ with the norm $\|[f]\| = \min_{c \in \mathbb{R}} \|f + c\|$. Without loss of generality, we always restrict to this quotient space and write $\|f\|$ to mean the norm in this quotient space. Moreover, for the quantities in equation (2.3) to exist, we shall restrict our attention to Banach spaces and require that differences in evaluations are continuous (i.e. the map $f \mapsto f(X_i) - f(X_j)$ is continuous for each $i$ and $j$). A Banach space is a normed vector space that is a complete metric space (see Ledoux and Talagrand (1991) and Royden (1988), chapter 10).

Borrowing terminology from game theory, we define two types of optimal designs: the *pure strategy optimal design (PSOD)* and the *mixed strategy optimal design (MSOD)*. For general $m \geqslant 2$, we define

$$B_{kk'}(W, f) = \frac{1}{p} \sum_{i:\, W_i = k} f(X_i) - \frac{1}{p} \sum_{i:\, W_i = k'} f(X_i).$$

The PSOD finds single assignments $W$ that on their own minimize these quantities.

*Definition 1.* Given subjects' baseline covariates $X \in \mathcal{X}^n$ and a magnitude function $\|\cdot\| : \mathcal{V} \to \mathbb{R} \cup \{\infty\}$, the *PSOD* chooses $W$ uniformly at random from the set of optimizers

$$W \in \arg\min_{W \in \mathcal{W}} \{M_{\mathrm{P}}^2(W) := \max_{\|f\| \leqslant 1} \max_{k \neq k'} B_{kk'}^2(W, f)\}.$$

We denote by $M_{\mathrm{P\text{-}opt}}^2$ the random variable equal to the optimal value.

The MSOD directly optimizes the distribution of assignments.

*Definition 2.* Given subjects' baseline covariates $X \in \mathcal{X}^n$ and a magnitude function $\|\cdot\| : \mathcal{V} \to \mathbb{R} \cup \{\infty\}$, the MSOD draws $W$ randomly according to a distribution $\sigma$ such that

$$\sigma \in \arg\min_{\sigma \in \Delta} \{M_{\mathrm{M}}^2(\sigma) := \max_{\|f\| \leqslant 1} \max_{k \neq k'} \sum_{W \in \mathcal{W}} \sigma(W) B_{kk'}^2(W, f)\}.$$

We denote by $M_{\mathrm{M\text{-}opt}}^2$ the random variable equal to the optimal value.

Both designs satisfy assumption 1. The PSOD does because of the symmetry of the objective function: if $W$ is optimal then any treatment permutation of $W$ is also optimal. The MSOD does by the construction of $\Delta$. This also implies that $M_{\text{M-opt}}^2 \leqslant M_{\text{P-opt}}^2$ (statewise dominance).

The objectives $M_{\text{P}}^2(W)$ and $M_{\text{M}}^2(\sigma)$ are the imbalance metrics that the designs seek to minimize. They are different in nature as one expresses imbalance of a particular assignment, in line with the common conception of the meaning of imbalance and, the other, the imbalance of a whole design. Since evaluation differences are linear and by assumption continuous, both $M_{\text{P}}^2(W)$ and $M_{\text{M}}^2(\sigma)$ are norms taken in the continuous dual Banach space and this guarantees that they are defined.

For mixed strategies, $M_{\text{M}}^2(\sigma)$ is actually determined by $n(n-1)/2$ sufficient statistics from $\sigma$.

*Theorem 2.* For any $\sigma \in \Delta$ define the matrix-valued $Q : \Delta \to \mathbb{R}^{n \times n}$ by

$$Q_{ij}(\sigma) = \sigma(\{W : W_i = W_j\}) - \frac{1}{m-1}\sigma(\{W : W_i \neq W_j\}).$$

Then, for any $k \neq k'$,

$$\sum_{W \in \mathcal{W}} \sigma(W) B_{kk'}^2(W, f) = \frac{2}{pn} \sum_{i,j=1}^n Q_{ij}(\sigma) f(X_i) f(X_j).$$

Consequently, $M_{\text{M}}^2(\sigma)$ depends on $\sigma$ only via $Q(\sigma)$:

$$M_{\text{M}}^2(\sigma) = M_{\text{M}}^2\{Q(\sigma)\} := \max_{\|f\| \leqslant 1} \frac{2}{pn} \sum_{i,j=1}^n Q_{ij}(\sigma) f(X_i) f(X_j).$$

For $m = 2$, feasible $Q$-matrices are $\mathcal{Q} := \text{convex-hull}(\mathcal{U})$ where $\mathcal{U} = \{u \in \{-1, 1\}^n : \Sigma_i u_i = 0\}$. All feasible matrices $Q \in \mathcal{Q}$ are positive semidefinite (symmetric with non-negative eigenvalues).

### 2.3. Structural information and existing designs and imbalance metrics

The definition of the optimal design above depends on a choice of norm $\|\cdot\|$. In this section, we take a reverse engineering approach to exhibit choices of these that recover various existing *a priori* balancing designs as optimal. In particular, this demonstrates the aptness and pervasiveness of this framework. In this section we consider two treatments, $m = 2$.

#### 2.3.1. Blocking and complete randomization

Randomized block designs are probably the most common non-completely randomized designs. In a complete-block design the sample is segmented into $b$ disjoint evenly sized blocks $\{i_{1,1}, \ldots, i_{1,2p_1}\}, \ldots, \{i_{b,1}, \ldots, i_{b,2p_b}\}$ so that baseline covariates are equal within each block and unequal between blocks ($X_{i_{l,j}} = X_{i_{l',j'}} \Leftrightarrow l = l'$). If any coarsening is done, we assume that it was done before and $X_i$ represents the coarsened value. The choice of coarsening is beyond the scope of this paper. Then complete randomization is applied to each block separately and independently of the other blocks. In incomplete blocks, there may be left-over subjects $i_{0,1}, \ldots, i_{0,b'}$ (e.g. a subject has no match) and one blocks subjects into evenly sized blocks so that there are fewest left-overs $b'$, breaking ties randomly about which subject is left over; complete randomization is then also applied to the left-overs.

Incomplete blocking maximizes a discrete measure of balance equal to the number of exact perfect matches across experimental groups. If complete blocking is feasible, then incomplete blocking necessarily recovers it. If all values of $X_i$ are distinct, then incomplete blocking is the same as complete randomization. As it is the most general, we shall treat only incomplete blocking. It turns out that incomplete blocking's exact matching metric corresponds to the space $L^\infty$.

*Theorem 3.* Let $\|f\| = \|f\|_\infty := \sup_{x \in \mathcal{X}} f(x)$. Then the PSOD is incomplete blocking.

This recovers as optimal both complete blocking (if feasible) and complete randomization (if $X_i \neq X_j$).

### 2.3.2. *Pairwise-matched allocation*

In pairwise-matched allocation, two treatments are considered, subjects are put into pairs to minimize the sum of pairwise distances in their covariates $\delta(X_i, X_j)$, and then each pair is split randomly between the two treatments (Greevy *et al.*, 2004). Usually $\delta$ is the pairwise Mahalanobis distance (for vector-valued covariates), but alternatives exist (see Barrios (2014), chapter 2). It turns out that the Lipschitz norm exactly recovers pairwise-matched allocation as an optimal design.

*Theorem 4.* Let a distance metric $\delta$ on $\mathcal{X}$ be given. Let

$$\|f\| = \|f\|_{\mathrm{lip}} = \sup_{x \neq x'} \left\{ f(x) - f(x') \right\} / \delta(x, x').$$

Then the PSOD is optimal pairwise-matched allocation with respect to the pairwise distance metric $\delta$.

Although $\|\cdot\|_{\mathrm{lip}}$ is only a seminorm on functions (i.e. $\|f\|_{\mathrm{lip}} = 0$ need not mean $f = 0$), in the quotient space with respect to constant functions, it is a norm and it forms a Banach space. Evaluation differences are well defined and continuous because $|f(X_i) - f(X_j)| \leqslant \|f\|_{\mathrm{lip}} \delta(X_i, X_j)$.

The motivation behind pairwise-matched allocation is that subjects with similar covariates should have similar outcomes. The interpretation of pairwise-matched allocation as a PSOD recasts this motivation as structure: that such a relationship holds between outcomes and covariates is the assumption of Lipschitz structure. Lipschitz continuous functions are essentially continuous functions of bounded variability and are extremely general. Unfortunately, they are so general that it is difficult to 'picture' the space of Lipschitz functions because it is not even separable (it has no countable dense subset). Comparing with blocking we see that, whereas blocking treats any two subjects with unequal covariates as potentially having expected outcomes that are as different as may be, pairwise-matched allocation presumes that unequal but similar covariates should lead to similar expected outcomes. This interpretation of pairwise-matched allocation also allows us to generalize it to $m \geqslant 3$ (see Section 4.1.2). Caliper allocation and an *a priori* version of the method of Kapelner and Krieger (2014) are also PSODs.

*Corollary 1.* Let $\delta_0 > 0$ and a distance metric $\delta$ be given. Define $\delta'(x, x') = \mathbb{I}_{[x \neq x']} \max \{\delta(x, x'), \delta_0\}$. Let $\|f\| = \|f\|_{\mathrm{lip}}$ with respect to $\delta'$. Then the PSOD is caliper allocation if it is feasible, i.e. choose at random from pairings that have all pairwise distances at most $\delta_0$, after blocking all exact matches.

*Theorem 5.* Let $\delta_0 > 0$ and a distance metric $\delta$ be given. Let $\|f\| = \max\{\|f\|_{\mathrm{lip}}, \|f\|_\infty / \delta_0\}$. Then the PSOD is as follows: minimize the sum of pairwise distances with respect to $\delta$ with the option of leaving a subject unmatched at a penalty of $\delta_0$ (thus no pairs further than $2\delta_0$ will be matched); matched pairs are randomly split between the two groups and unmatched subjects are completely randomized.

This method mimics optimal pairwise-matched allocation for subjects who are close to one another but avoids pairing subjects who would be badly matched for whom there is little to suggest a similar response to treatments and it may be more robust to randomize these completely. Taking $\delta_0 \to \infty$ recovers pariwise-matched allocation whereas taking $\delta_0 \to 0$ recovers complete randomization.

### 2.3.3. Rerandomization

The method of Morgan and Rubin (2012) formalizes the common, but arguably historically haphazard, practice of rerandomization as a principled, theoretically grounded *a priori* balancing method. Morgan and Rubin (2012) considered two treatments: vector-valued baseline covariates $\mathcal{X} = \mathbb{R}^d$, and an imbalance metric equal to the Mahalanobis distance between group means,

$$M_{\text{Rerand}}^2(W) = \left( \frac{2}{n} \sum_{i=1}^n u_i X_i \right)^{\text{T}} \hat{\Sigma}^{-1} \left( \frac{2}{n} \sum_{i=1}^n u_i X_i \right), \tag{2.4}$$

where $\hat{\Sigma}$ is the sample covariance matrix of $X$. They reinterpreted rerandomization as a heuristic algorithm that repeatedly draws random $W$ to solve the constraint satisfaction problem $\exists? W : M_{\text{Rerand}}^2 \leqslant t$ for a given $t$. They also proposed a normal approximation method for selecting $t$ to correspond to a particular acceptance probability of a random $W$.

We can recover equation (2.4) by using our framework. Let $\mathcal{F} = \text{span}\{1, x_1, \ldots, x_d\}$ and define $\|f\|^2 = \beta^{\text{T}} \hat{\Sigma} \beta + \beta_0^2$ for $f(x) = \beta_0 + \beta^{\text{T}} x$. Using duality of norms,

$$M_{\text{P}}^2(W) = \max_{\|f\| \leqslant 1} B^2(W, f) = \left\{ \max_{\beta^{\text{T}} \hat{\Sigma} \beta \leqslant 1} \beta^{\text{T}} \left( \frac{n}{2} \sum_{i=1}^n u_i X_i \right) \right\}^2 = M_{\text{Rerand}}^2(W).$$

Alternatively, we can let $\|f\|^2 = \beta^{\text{T}} \beta + \beta_0^2$ and normalize the data in advance.

Morgan and Rubin (2012) argued that, when a linear model is known to hold, i.e.

$$Y_{ik} = \beta_0 + \beta^{\text{T}} X_i + \tau \mathbb{I}(k=1) + \epsilon_i \qquad i = 1, \ldots, n, \quad k = 1, 2, \tag{2.5}$$

then fixing $t$ and rerandomizing until $M_{\text{Rerand}}^2(W) \leqslant t$ yields a reduction in variance relative to complete randomization that is constant over $n$ (for large $n$):

$$1 - \text{var}(\hat{\tau})/\text{var}(\hat{\tau}^{\text{CR}}) \approx \eta\{1 - \text{var}(\epsilon_i)/\text{var}(Y_{i1})\}, \qquad \eta \in (0, 1) \text{ constant over } n.$$

For us, the imbalance metric is a direct consequence of structure (equation (2.5) implies $f_k \in \mathcal{F}$) and *fully* minimizing it leads to the nearly best conceivable reduction in variance for moderate $n$ (corollary 2, Section 3.3):

$$1 - \text{var}(\hat{\tau})/\text{var}(\hat{\tau}^{\text{CR}}) \to 1 - \text{var}(\epsilon_i)/\text{var}(Y_{i1}) \qquad \text{at a linear rate } 2^{-\Omega(n)}.$$

It is important to keep in mind, however, that the assumption that such a finite dimensional linear model (2.5) is valid is a parametric, and therefore fragile, assumption (see example 2 in Section 2.4.1).

### 2.3.4. Other finite dimensional spaces and the method of Bertsimas et al. (2015)

Consider any norm on a finite dimensional subspace $\mathcal{F}$, which can be written as $\mathcal{F} = \text{span}\{\phi_1, \ldots, \phi_r\}$ for some $\phi_i : \mathcal{X} \to \mathbb{R}$. Any such space is always a Banach space and evaluations are always continuous (Hunter and Nachtergaele (2001), theorems 5.33 and 5.35). An important example is the $q$-norm. The $q$-norm on $\mathbb{R}^r$ is $\|\beta\|_q = (\Sigma_i |\beta_i|^q)^{1/q}$ for $1 \leqslant q < \infty$ and $\|\beta\|_\infty = \max_i |\beta_i|$. We can use this to define a norm on $\mathcal{F}$: for $f = \beta_1 \phi_1 + \ldots + \beta_r \phi_r \in \mathcal{F}$, we let $\|f\| = \|\beta\|_q$. This yields

$$M_{\mathrm{P}}^2(W) = \left\| \left( \frac{n}{2} \sum_{i=1}^{n} u_i \phi_1(X_i), \ldots, \frac{n}{2} \sum_{i=1}^{n} u_i \phi_r(X_i) \right) \right\|_{q^*}^2$$

where $q^*$ satisfies $1/q + 1/q^* = 1$. Hence, the optimal design matches the sample $\phi_j$-moments between the groups by minimizing a norm in the vector of the $r$-moment mismatches. In Section 2.3.3, we saw that a scaled 2-norm on $\mathcal{F} = \mathrm{span}\{1, x_1, \ldots, x_d\}$ gave rise to a groupwise Mahalanobis metric. Fixing $\rho > 0$ and endowing $\mathcal{F} = \mathrm{span}\{1, x_1, \ldots, x_d, x_1^2/\rho, \ldots, x_d^2/\rho, x_1 x_2/(2\rho), \ldots, x_{d-1} x_d/(2\rho)\}$ with the $\infty$-norm and normalizing the data will recover the method of Bertsimas *et al.* (2015), which optimizes the balance in covariate means and centred second moments by using mixed integer programming.

### 2.4. Kernel allocation

In the previous section, we saw how various well-known imbalance metrics and designs were recovered as optimal under certain norms. In this section, we consider norms that are induced by kernels and study the new class of optimal designs that arise from these, which we call kernel allocation. We explore the implications on practice in Section 6. In this section, we treat general $m \geqslant 2$.

We shall express structure by using reproducing kernel Hilbert spaces (RKHSs). A Hilbert space is an inner product space such that the norm that is induced by the inner product, $\|f\|^2 = \langle f, f \rangle$, yields a Banach space. An RKHS $\mathcal{F}$ is a Hilbert space of functions for which, for every $x \in \mathcal{X}$, the map $f \mapsto f(x)$ from functions to their value at $x$ is a continuous mapping (Berlinet and Thomas-Agnan, 2004). Continuity and the Riesz representation theorem imply that for each $x \in \mathcal{X}$ there is $\mathcal{K}(x, \cdot) \in \mathcal{F}$ such that $\langle \mathcal{K}(x, \cdot) f(\cdot) \rangle = f(x)$ for every $f \in \mathcal{F}$. The symmetric map $\mathcal{K} : \mathcal{X} \times \mathcal{X} \to \mathbb{R}$ is called the reproducing kernel of $\mathcal{F}$. The name is motivated by the fact that $\mathcal{F} = \mathrm{closure} \, [\mathrm{span}\{\mathcal{K}(x, \cdot) : x \in \mathcal{X}\}]$. Thus $\mathcal{K}$ fully characterizes $\mathcal{F}$. Prominent examples of kernels for $\mathcal{X} \subset \mathbb{R}^d$ are as follows.

(a) The linear kernel $\mathcal{K}(x, x') = x^{\mathrm{T}} x'$: the RKHS spans the finite dimensional space of linear functions and induces a 2-norm on coefficients.
(b) The polynomial kernel $\mathcal{K}_s(x, x') = (1 + x^{\mathrm{T}} x'/s)^s$: the RKHS spans the finite dimensional space of all polynomials of degree up to $s$.
(c) Any kernel $\mathcal{K}(x, x') = \Sigma_{i=0}^{\infty} a_i (x^{\mathrm{T}} x')^i$ with $a_i \geqslant 0$ (subject to convergence): this includes the previous two examples. Another case is the exponential kernel $\mathcal{K}(x, x') = \exp(x^{\mathrm{T}} x')$, which can be seen as the infinite dimensional limit of the polynomial kernel. The corresponding RKHS is infinite dimensional (non-parametric).
(d) The Gaussian kernel $\mathcal{K}(x, x') = \exp(-\|x - x'\|^2)$: the corresponding RKHS is infinite dimensional (non-parametric) and was studied in Steinwart *et al.* (2006).

For given $X \in \mathcal{X}^n$ and an RKHS with kernel $\mathcal{K}$, we shall often use the Gram matrix $K_{ij} = \mathcal{K}(X_i, X_j)$. The Gram matrix is always positive semidefinite and as such it has a matrix square root $K = \sqrt{K}\sqrt{K}$.

Some infinite dimensional RKHSs, known as universal, can arbitrarily approximate any continuous function, making them incredibly general. The Gaussian kernel is universal and the exponential kernel is universal on compact spaces (Sriperumbudur *et al.*, 2011).

*Definition 3.* For $\mathcal{X}$ Hausdorff, an RKHS $\mathcal{F}$ (or the corresponding kernel) is $C_0$ *universal* if, for any continuous function $g : \mathcal{X} \to \mathbb{R}$ with compact support (i.e. $\{x : g(x) \neq 0\} \subset C$ for $C$ compact) and $\eta > 0$, there exists $f \in \mathcal{F}$ such that $\sup_{x \in \mathcal{X}} |f(x) - g(x)| \leqslant \eta$.

An RKHS is a class of functions. Therefore, in our framework, any kernel induces imbalance metrics and optimal designs. We call these designs kernel allocation.

**Theorem 6.** Let $\mathcal{F}$ be an RKHS with kernel $\mathcal{K}$. Then,

$$M_{\mathrm{P}}^2(W) = \frac{1}{p^2} \max_{k \neq k'} \sum_{i,j=1}^{n} (w_{ik} - w_{ik'}) K_{ij} (w_{jk} - w_{jk'}), \tag{2.6}$$

and

$$M_{\mathrm{M}}^2(Q) = \frac{2}{np} \lambda_{\max}(\sqrt{K} Q \sqrt{K}). \tag{2.7}$$

The problem of minimizing expressions (2.6) or (2.7) can be interpreted as a multiway multicriterion number partitioning problem. In particular, for the case when $m = 2$, $\mathcal{X} = \mathbb{R}$ and $\mathcal{K}(x, x') = xx'$ ($K = XX^{\mathsf{T}}$), we recover the standard balanced number partitioning problem. Recalling that $\mathcal{Q} = \text{convex-hull}(\mathcal{U})$, we see that

$$\begin{aligned}
\frac{n}{2} M_{\mathrm{P\text{-}opt}} &= \sqrt{\{\min_{u \in \mathcal{U}} u^{\mathsf{T}}(XX^{\mathsf{T}})u\}} = \min_{u \in \mathcal{U}} \left| \sum_{i=1}^{n} u_i X_i \right|, \\
\frac{n}{2} M_{\mathrm{M\text{-}opt}} &= \sqrt{[\min_{Q \in \mathcal{Q}} \text{tr}\{Q(XX^{\mathsf{T}})\}]} = \min_{u \in \mathcal{U}} \left| \sum_{i=1}^{n} u_i X_i \right|,
\end{aligned} \tag{2.8}$$

where the last equality is due to the facts that $\lambda_{\max}(M) = \text{tr}(M)$ if $M$ is rank 1 symmetric and that a linear objective on a polytope is optimized at a corner point. This reduction shows that both problems are 'NP hard' (see Garey and Johnson (1979), problem [SP12]). This also shows that the MSOD is only relevant for high rank $K$. Universal kernels always have full rank $K$ for distinct $X$.

Such partitioning problems generically have unique optima up to permutation so the PSOD usually randomizes among the $m!$ permutations of a single partition of subjects. This is not generally so for the MSOD. Consider $m = 2$. Since the affine hull of $\mathcal{U}$ is $n - 1$ dimensional, the MSOD mixes at least rank$(K) - 1$ distinct partitions (and permutations thereof). Moreover, by Carathéodory's theorem any $Q \in \mathcal{Q}$ can be identified as the convex combination of $n(n-1)$ points in $\{uu^{\mathsf{T}} : u \in \mathcal{U}\}$ (whose affine hull is $n(n-1) - 1$ dimensional) so that the MSOD objective $M_{\mathrm{M}}^2(\sigma)$ of any *a priori* balancing design $\sigma \in \Delta$ can also be achieved by mixing no more than $n(n-1)$ distinct partitions.

In Sections 4.1.3 and 4.2 we shall study how we solve the PSOD and MSOD respectively. For now let us consider two concrete examples with the various designs we have so far studied.

### 2.4.1. Example 2
Consider the following set-up: we measure $d \geqslant 2$ baseline covariates for each subject that are uniformly distributed in the population $X_i \sim \text{Unif}([-1, 1]^d)$, the two treatments ($m = 2$) have constant individual effects $Y_{i1} - Y_{i2} = \tau$, and the conditional expectation of outcomes depends on two covariates only: $\mathbb{E}[Y_{i1} | X = x] - \tau/2 = \mathbb{E}[Y_{i2} | X = x] + \tau/2 = \hat{f}(x) = \hat{f}(x_1, x_2)$. We consider a variety of conditional expectation functions:

(a) linear, $\hat{f}(x_1, x_2) = x_1 - x_2$;
(b) quadratic, $\hat{f}(x_1, x_2) = x_1 - x_2 + x_1^2 + x_2^2 - 2x_1 x_2$;
(c) cubic, $\hat{f}(x_1, x_2) = x_1 - x_2 + x_1^2 + x_2^2 - 2x_1 x_2 + x_1^3 - x_2^3 - 3x_1^2 x_2 + 3x_1 x_2^2$;
(d) sinusoidal, $\hat{f}(x_1, x_2) = \sin(\pi/3 + \pi x_1/3 - 2\pi x_2/3) - 6\sin(\pi x_1/3 + \pi x_2/4) + 6\sin(\pi x_1/3 + \pi x_2/6)$.

We do not consider the case of no relationship ($\hat{f}(x_1, x_2) = c$) because theorem 9 in Section 3.1 proves that in this case any *a priori* balancing design yields the exact same estimation variance.

To simulate the common situation where some covariates matter and some do not, and which is not known *a priori*, we consider both the case $d = 2$ (only balance the relevant covariates) and $d = 4$ (also balance some covariates that turn out to be irrelevant).

We consider the following designs:

(a) complete randomization, i.e. the PSOD for $L^\infty$;
(b) blocking on the orthant of $X_i$ ($d$ two-level factors), i.e. the PSOD for $L^\infty$ after coarsening;
(c) rerandomization with 1% (exact) acceptance probability and Mahalanobis objective;
(d) pairwise-matched allocation with Mahalanobis distance, i.e. the PSOD for the Lipschitz norm;
(e) linear kernel allocation PSOD;
(f) quadratic kernel allocation PSOD (polynomial kernel with $s = 2$);
(g) Gaussian kernel allocation MSOD;
(h) exponential kernel allocation MSOD (for the MSODs we use the heuristic solution given by algorithm 3 in Section 4).

All these designs result in an unbiased estimate of $\mathrm{SATE} = \mathrm{PATE} = \tau$ and can therefore be compared on their variance. In Fig. 2 we plot the variances of the resulting estimators relative to $V_n = \mathrm{var}(\mathrm{SATE}) + \mathrm{var}(\epsilon_{11} + \epsilon_{12})/n$.

There are several features to note. One is that, when a parametric model is correctly specified and specifically optimized for, the variance (relative to $V_n$) shrinks linearly (inverse exponentially)—we argue that this is a general phenomenon in Section 3.3. This phenomenon is clearest in the case of linear kernel allocation under a linear conditional expectation. But, linear kernel allocation does not do so well when the linear model is misspecified. Quadratic kernel allocation also has a fast linear rate of convergence under the linear conditional expectation, but it performs much better in all the other cases, both when a quadratic model is correctly specified and when it is not. Gaussian and exponential kernel allocation seem to have uniformly good performance in all cases and in particular still exhibit what would seem to be linear convergence for the linear and quadratic cases. It would seem that these non-parametric methods strike a good compromise between efficiency and robustness. Finally, we note that, compared with balancing only relevant covariates ($d = 2$), balancing also irrelevant covariates ($d = 4$) leads to some loss of efficiency since it leads to a less good match for the relevant covariates, but the order of the rate of convergence (linear) is the same.

### 2.4.2. Example 3

We now consider the effect of *a priori* balance on a real data set. We use the diabetes study data set from Efron *et al.* (2004) described therein as follows:

'Ten [$d = 10$] baseline variables [denoted $X_i$], age, sex, body mass index, average blood pressure, and six blood serum measurements were obtained for each of [$n = 442$] diabetes patients, as well as the response of interest [denoted $Y_i$], a quantitative measure of disease progression one year after baseline'.

We consider a hypothetical experiment where the prognostic features $X_i$ are measured at the outset, a control or treatment is applied, and the response after 1 year is measured. In our hypothetical set-up, the treatment reduces disease progression by exactly $\tau$ so that $Y_{i1} = Y_i'$ and $Y_{i2} = Y_i' - \tau$. Fixing $n$, we draw $n$ subjects with replacement from the population of 442, normalize the covariate data so that the sample of $n$ has zero sample mean and identity sample

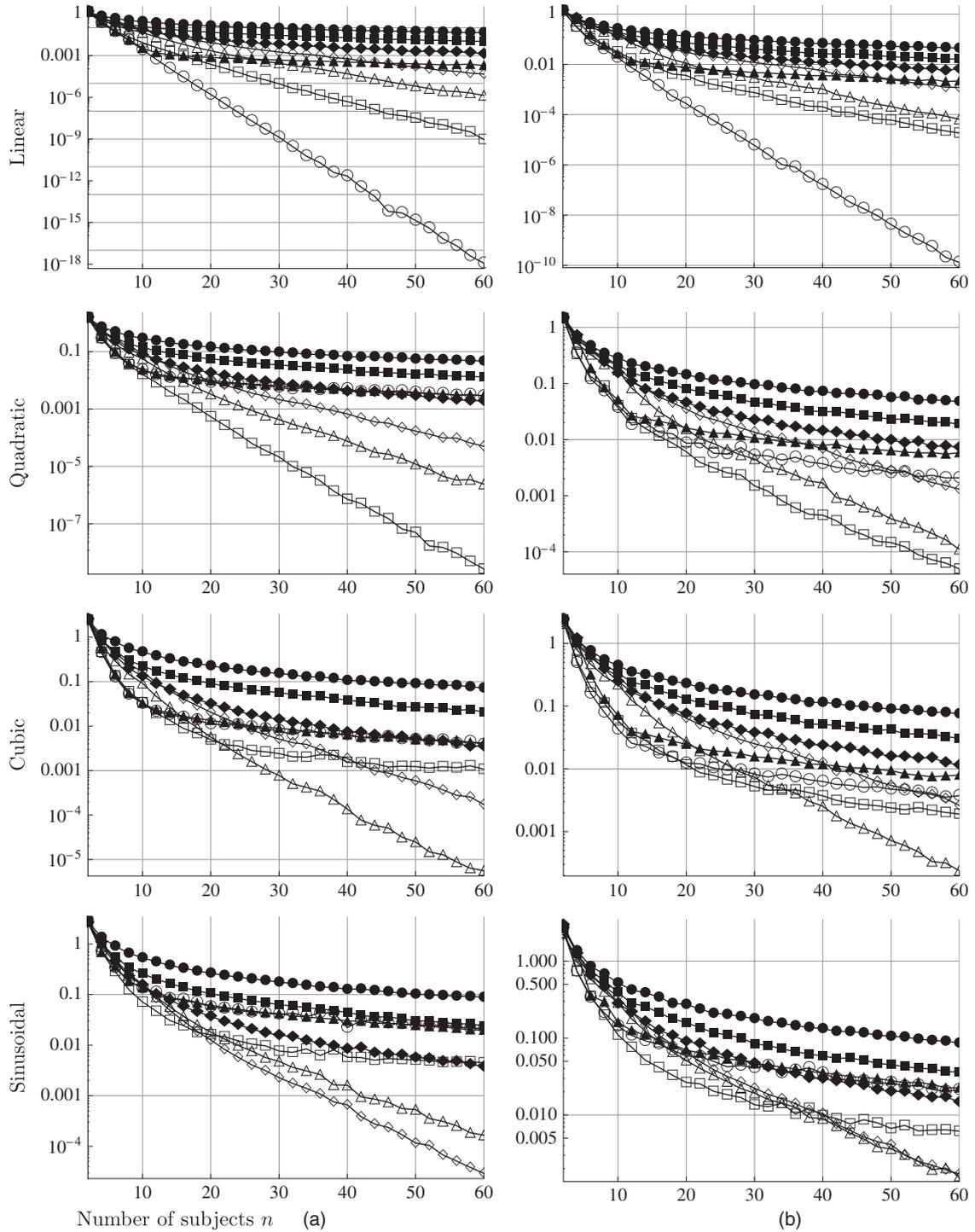

**Fig. 2.** Estimation variance $\mathrm{var}(\hat{\tau}) - V_n$ under various designs, covariate dimensions and conditional expectations in example 2 (●, optimal $L^\infty$ (complete randomization); ○, optimal linear; ■, optimal $L^\infty$ coarsened (blocking); □, optimal quadratic; ◆, optimal Lipschitz (pairwise match); ◇, optimal Gaussian; ▲, rerandomization; △, optimal exponential): (a) $d = 2$; (b) $d = 4$

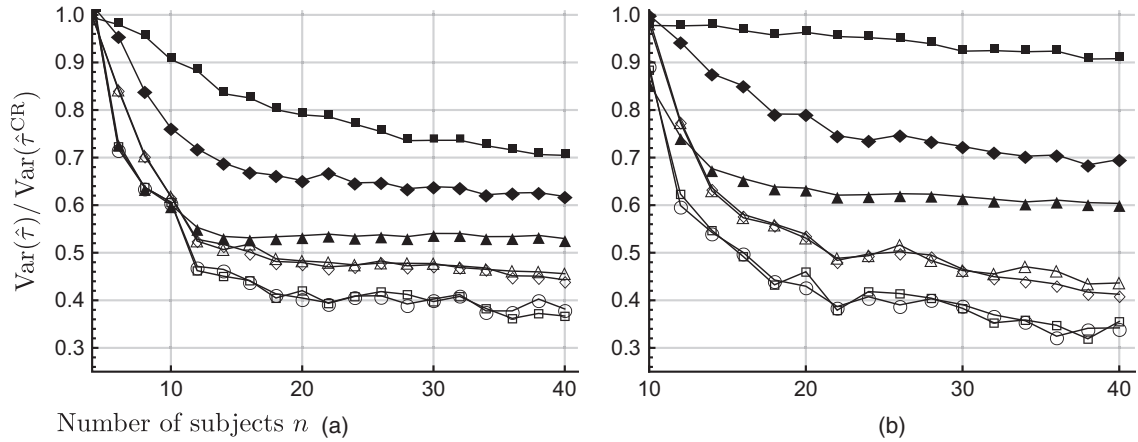

**Fig. 3.** Relative variance var($\hat{\tau}$)/var($\hat{\tau}^{CR}$) under various designs for the diabetes data set in example 3 (■, optimal $L^\infty$ coarsened (blocking); ○, optimal linear; ♦, optimal Lipschitz (pairwise matching); □, optimal quadratic; ▲, rerandomization; ◇, optimal Gaussian; △, optimal exponential): (a) $d = 4$; (b) $d = 10$

covariance and divide by $d = 10$, apply each of the *a priori* balancing designs that were considered in example 2 to the normalized covariates, and finally apply the treatments and measure the responses and the mean differences $\hat{\tau}$. Again, we consider either balancing all $d = 10$ covariates or only the $d = 4$ covariates that were ranked most important by Efron *et al.* (2004): $\{3, 9, 4, 7\}$. We plot estimation variances relative to complete randomization in Fig. 3.

For larger $n$, the relative variance of each method stabilizes around a particular ratio. Each of blocking, pairwise-matched allocation and rerandomization result in a higher ratio when attempting to balance all covariates compared with balancing only the four most important. For example, rerandomization on all 10 covariates gives about 60% of complete randomization's variance whereas restricting to the important covariates yields about 53%. These classic balancing designs appear to do worse when given additional covariates instead of exploiting the additional prognostic content that is available.

In contrast, kernel allocation yields lower relative variances for both $d = 10$ and $d = 4$, converging slower for $d = 10$ but using the small additional prognostic content of the extra covariates to reduce the variance even further. For example, linear and quadratic kernel allocation both yield about 40% of complete randomization's variance for $d = 4$ and about 35% for $d = 10$, taking only slightly longer to fall below about 40% when $d = 10$. This can be attributed to the linear rate at which the optimal designs eliminate imbalances (see Section 3.3). Thus, even if there are some less relevant variables, all are quickly nearly perfectly balanced for modest $n$; the only limiting factors are the residuals $\epsilon_{ik}$, which, by definition, cannot be controlled for by using the covariates $X$ alone (see corollary 3 in Section 3.1).

### 2.4.3. *Aside: worst-case analysis, randomization, pure strategy optimal design versus mixed strategy optimal design and a Bayesian interpretation*

In the preceding sections, we motivated a variety of designs as guarding against the worst possible realization by minimizing the worst-case conditional variance. A question that arises is, why a worst-case analysis? Fisher (1990) once justified randomization by considering a game with the Devil (or, nature) where randomization was a means to maintain unbiasedness and to attain efficient estimates in spite of the Devil's machinations (Senn, 1994). The worst-case point of view that was presented in Sections 2.1 and 2.2 to motivate the designs is in line with these views of Fisher. As we saw in Section 2.3, a worst-case point of view also reproduces many standard designs and imbalance measures used in practice.

The aspect of randomization that achieves unbiasedness in spite of uncertainty, known as blinding, is the gist of assumption 1, part (b) (see Senn (1994)). It may be appealing to think that, to reduce variance the most, while blinding the identity of treatments, one should pick the single most balanced partition of subjects (which is then permuted randomly over treatments for blinding). This, however, is not necessarily right from a worst-case point of view. We saw in theorem 1 and example 1 that, in the absence of structure, taking rerandomization to the extreme of minimization led to the worst possible variance and only a design that takes randomization to the extreme of complete randomization was optimal. This is true even when (certain) structure is present: whereas randomizing treatment identities eliminates bias in the face of uncertainty, randomizing over partitions can reduce variance in the face of uncertainty. If the objective is to minimize variance in the face of uncertainty restricted by structure, the correct objective is $M_M^2(\sigma)$. In Section 2.4, we argued that if $\mathcal{F}$ is an RKHS and $\mathrm{rank}(K) > 1$ then the MSOD must randomize between more than one partition (more than two assignments for $m = 2$). (In contrast, if $\mathrm{rank}(K) = 1$ then equation (2.8) shows that the MSOD randomizes only a single partition, generically.) Nonetheless, choosing single partitions that on their own have good balance, rather than mixtures of partitions that achieve minimax, has practical appeal, especially in interpretability, and can achieve similar performance in practice as we saw in example 3.

In particular, the kernel allocation PSOD, which generically chooses only one partition, can be interpreted as a Bayes optimal design. The interpretation is similar to the Bayesian interpretation of ridge regression (see Kimeldorf and Wahba (1970) and Rasmussen and Williams (2006) section 6.2). Let $m = 2$ and let $\mathcal{F}$ be a given RKHS with kernel $\mathcal{K}$. Let us assume a Gaussian prior on $\hat{f}$ with covariance operator $\mathcal{K}$, i.e. $\hat{f}(x)$ is Gaussian for every $x \in \mathcal{X}$ and the covariance of $\hat{f}(x)$ and $\hat{f}(x')$ is equal to $\mathcal{K}(x, x')$. Then, using equation (2.2), we have that the Bayes variance risk of $\hat{\tau}$ is

$$\mathbb{E}[B^2(W, f)|X, Y] = \frac{4}{n^2} \sum_{i,j=1}^{n} u_i u_j \mathbb{E}[f(X_i) f(X_j)|X, Y] = \frac{4}{n^2} \sum_{i,j=1}^{n} u_i u_j \mathcal{K}(X_i, X_j) = M_P^2(W).$$

So, the kernel allocation PSOD minimizes Bayes risk given the data. Note that randomization is unnecessary from a standard Bayesian perspective (Kadane and Seidenfeld, 1990; Savage, 1964). So, *not* blinding assignment (not permuting the optimal partition) is *also* Bayes optimal but may be biased.

### 2.4.4. Aside: agnostic notes on regression adjustments and kernel allocation

In estimating PATE, the assumption that $\mathcal{F}$ is a finite dimensional (parametric) RKHS enables a regression adjustment whereby we estimate the parametric regression model to obtain the population effect of treatment assignment. Freedman (1997, 2008), Robins (1994, 2004) and Berk (2004) have drawn attention to the fact that, when such a parametric modelling assumption is false, this estimate may be extremely inaccurate and may lead to false conclusions. In contrast, Lin (2013) pointed out that a careful adjustment that includes interaction terms (equivalently, fitting a regression separately to the treated and control samples) will incur only slight bias and will be asymptotically valid.

In contrast, if we use kernel allocation, then, *regardless* of the validity of the parametric RKHS model and for any sized sample, our simple estimator $\hat{\tau}_{kk'}$ is still completely unbiased for SATE in finite (fixed) samples and for PATE in repeated draws and testing inference is valid. In the case that a parametric RKHS is correctly specified, we show in Section 3.3 that any error in $\hat{\tau}_{kk'}$ that is explainable by $X$ (i.e. the error that is adjustable by regression) will vanish at a linear rate $2^{-\Omega(n)}$, achieving the same efficient variance of adjustment without incurring bias. In this way we have the best of both worlds: no bias or invalid inference *regardless* of sample size or valid specification and nearly perfect control of explained variance for even quite small samples under a valid parametric specification, such as linear.

## 3. Characterizations of *a priori* balancing designs

In this section, we theoretically characterize the variance, consistency and convergence rates of estimation under kernel allocation and other optimal designs.

### 3.1. Variance

We begin by decomposing the variance of the estimator under any *a priori* balancing design.

*Theorem 7.* Suppose that assumption 1 is satisfied. Then, for all $k \neq k'$:

(a) $\hat{\tau}_{kk'}$ is conditionally and marginally unbiased, $\mathbb{E}[\hat{\tau}_{kk'}|X, Y] = \text{SATE}_{kk'}$ and $\mathbb{E}[\hat{\tau}_{kk'}] = \text{PATE}_{kk'}$;

(b) $\hat{\tau}_{kk'} = \text{SATE}_{kk'} + D_{kk'} + E_{kk'}$, where

$$D_{kk'} := \frac{1}{m}\sum_{l \neq k} B_{kl}(W, f_k) - \frac{1}{m}\sum_{l \neq k} B_{k'l}(W, f_{k'}),$$

$$E_{kk'} := \frac{1}{n}\sum_{i=1}^{n}\{(mw_{ik} - 1)\epsilon_{ik} - (mw_{ik'} - 1)\epsilon_{ik'}\};$$

(c) $\text{SATE}_{kk'}$, $D_{kk'}$ and $E_{kk'}$ are all uncorrelated so that

$$\text{var}(\hat{\tau}_{kk'}) = \frac{1}{n}\text{var}(Y_{1k} - Y_{1k'}) + \text{var}(D_{kk'}) + \frac{1}{n}\text{var}(\epsilon_{1k} + \epsilon_{1k'}) + \frac{m-2}{n}\{\text{var}(\epsilon_{1k}) + \text{var}(\epsilon_{1k'})\}. \tag{3.1}$$

In equation (3.1), every term except $\text{var}(D_{kk'})$ is unaffected by our choice of *a priori* balancing design. This can be seen because, in equation (3.1), $W$ appears nowhere except in $D_{kk'}$. Below we provide a bound on $\text{var}(D_{kk'})$ based on the expected minimal imbalance.

*Theorem 8.* For any $\mathcal{F}$, if the PSOD or MSOD is used,

$$\text{var}(D_{kk'}) \leqslant \frac{(\|f_k\| + \|f_{k'}\|^2)}{2}\left(1 - \frac{1}{m}\right)\mathbb{E}[M_{\text{opt}}^2] \propto \mathbb{E}[M_{\text{opt}}^2], \tag{3.2}$$

where $M_{\text{opt}}^2 = M_{\text{P-opt}}^2$ or $M_{\text{opt}}^2 = M_{\text{M-opt}}^2$ respectively.

In inequality (3.2), $(\|f_k\| + \|f_{k'}\|)^2$ is unknown but constant, merely scaling the bound. Combining these intimately connects imbalance before treatment and randomization to estimation variance afterwards.

*Corollary 2.* For any $\mathcal{F}$, if the PSOD or MSOD is used,

$$\text{var}(\hat{\tau}_{kk'}) \leqslant \frac{1}{n}\text{var}(Y_{1k} - Y_{1k'}) + \frac{(\|f_k\| + \|f_{k'}\|)^2}{2}\left(1 - \frac{1}{m}\right)\mathbb{E}[M_{\text{opt}}^2]$$
$$+ \frac{1}{n}\text{var}(\epsilon_{1k} + \epsilon_{1k'}) + \frac{m-2}{n}\{\text{var}(\epsilon_{1k}) + \text{var}(\epsilon_{1k'})\}.$$

*Corollary 3.* Suppose that $m = 2$ and that individual effects are constant, $Y_{i1} - Y_{i2} = \tau$. Let $\sigma_{\text{total}}^2 = \text{var}(Y_{i1}) = \text{var}(Y_{i2})$, $\sigma_{\text{residual}}^2 = \text{var}(\epsilon_{i1}) = \text{var}(\epsilon_{i2})$ and $R^2 = 1 - \sigma_{\text{residual}}^2 / \sigma_{\text{total}}^2$ (the fraction of $Y$-variance explained by $X$). The variance under the PSOD or MSOD relative to the variance under complete randomization satisfies

$$1 - R^2 \leqslant \frac{\text{var}(\hat{\tau})}{\text{var}(\hat{\tau}^{\text{CR}})} \leqslant 1 - R^2 + \frac{n}{16\sigma_{\text{total}}^2}(\|f_1\| + \|f_2\|)^2 \mathbb{E}[M_{\text{opt}}^2].$$

Alternatively, the relative *reduction* in variance is 1 minus this inequality above. Despite the constant effect assumption, this bound provides important insights. On the one hand, it says that any *a priori* balancing effort can never do better than $1 - R^2$ relative to no balancing (complete randomization). This makes sense: balancing based on $X$ alone can help only to the extent that $X$ is predictive of outcomes. On the other hand, it says that, if $\mathbb{E}[M_{\text{opt}}^2]$ decays superlogarithmically, i.e. $o(1/n)$, then the relative variance converges to the best theoretically possible, which is $1 - R^2$. In Section 3.3 we study a case where the convergence is linear, i.e. $2^{-\Omega(n)}$, which is much faster than logarithmic.

Another conclusion that we can draw from theorem 7 is that balancing irrelevant baseline covariates leads to no loss of efficiency regardless of any *a priori* balancing—any design has the same variance. This is a complete generalization of special results on no loss of efficiency in blocking (Cochran and Cox (1957), section 4.26) or pairwise-matched allocation (Chase, 1968) on irrelevant covariates. (Note, however, that excising irrelevant covariates may improve balance in relevant covariates; see also example 3.)

*Theorem 9.* Suppose that assumption 1 is satisfied and that $Y_{1k}$ and $Y_{1k'}$ are mean independent of $X_1$, i.e. $\mathbb{E}[Y_{1k}|X_1] = \mathbb{E}[Y_{1k}]$ and $\mathbb{E}[Y_{1k'}|X_1] = \mathbb{E}[Y_{1k'}]$. Then, $\text{var}(\hat{\tau}_{kk'}) = \text{var}(\hat{\tau}_{kk'}^{\text{CR}})$.

Finally, we treat the case where the baseline covariates are relevant but we have ruled out the true conditional expectation function, i.e. $f_k \notin \mathcal{F}$. When $f_k \notin \mathcal{F}$ we have $\|f_k\| = \infty$ and the bound (3.2) is trivial. Accounting for the distance between $f_k$ and $\mathcal{F}$, an alternative bound is possible.

*Theorem 10.* For any $\mathcal{F}$, if the PSOD or MSOD is used,

$$\text{var}(D_{kk'}) \leqslant \left(1 - \frac{1}{m}\right) \inf_{g_k, g_{k'} \in \mathcal{F}} \left\{ (\|g_k\| + \|g_{k'}\|)^2 \mathbb{E}[M_{\text{opt}}^2] + \frac{2}{m}(\|f_k - g_k\|_2 + \|f_{k'} - g_{k'}\|_2)^2 \right\},$$

where $\|g\|_2^2 = \mathbb{E}[g(X_1)^2]$ is the $L^2$-norm in the measure of $X_1$.

Note that the above bound is finite since, by the assumption that potential outcomes have second moments, we have $\|f_k\|_2 < \infty$. A universal RKHS is an example of a space that can approximate *any* function arbitrarily closely in $L^2$, yielding the following theorem as a corollary.

*Theorem 11.* For any $C_0$-universal RKHS $\mathcal{F}$, if the PSOD or MSOD is used and $\mathcal{X}$ is locally compact (e.g. $\mathbb{R}^d$), then for any $\eta > 0$ there are $g_k, g_{k'} \in \mathcal{F}$ such that for all $n$

$$\text{var}(D_{kk'}) \leqslant \left(1 - \frac{1}{m}\right)(\|g_k\| + \|g_{k'}\|)^2 \mathbb{E}[M_{\text{opt}}^2] + \eta. \tag{3.3}$$

Note that we do not require $f_k$ and $f_{k'}$ to be continuous and that the choice of $g_k$ and $g_{k'}$ need not change with $n$. So, as $\mathbb{E}[M_{\text{opt}}^2] \to 0$ with $n \to \infty$, only $\eta$ remains in inequality (3.3) and $\eta$ is arbitrary (see theorem 13 in Section 3.2).

### 3.2. Consistency

An estimator is said to be consistent if it converges to the quantity that it tries to estimate. Using laws of large numbers in Banach spaces, we can express conditions for consistency of the PSOD and MSOD. (Note that, since $\text{SATE}_{kk'} \to \text{PATE}_{kk'}$ almost surely, consistency in estimating SATE or PATE is the same.)

*Definition 4.* A Banach space is said to be *B convex* if there exists $N \in \mathbb{N}$ and $\eta < N$ such that for every $g_1, \ldots, g_N$ with $\|g_i\| \leqslant 1 \, \forall i$ there is a choice of signs so that $\|\pm g_1 \pm \ldots \pm g_N\| \leqslant \eta$.

All the Banach spaces so far considered have been *B* convex with the exception of $L^\infty$. Generally, any Hilbert space or finite dimensional Banach space is *B* convex (Ledoux and Talagrand (1991), chapter 9).

*Theorem 12.* Suppose that $f_k, f_{k'} \in \mathcal{F}$. Then, under either the PSOD or MSOD, we have that $\hat{\tau}_{kk'} - \text{SATE}_{kk'} \to 0$ almost surely if either of the following conditions hold:

(a) $\mathcal{F}$ is *B* convex and $\mathbb{E}[\max_{\|f\| \leqslant 1} \{f(X_1) - f(X_2)\}^2] < \infty$ or
(b) $\mathcal{F}$ is a Hilbert space and $\mathbb{E}[\max_{\|f\| \leqslant 1} \{f(X_1) - f(X_2)\}] < \infty$.

If we use a universal RKHS, we can guarantee consistency regardless of misspecification $f_k$, $f_{k'} \notin \mathcal{F}$.

*Theorem 13.* Suppose that $\mathcal{F}$ is a $C_0$-universal RKHS and $\mathbb{E}[\max_{\|f\| \leqslant 1} \{f(X_1) - f(X_2)\}] < \infty$. Then, under either the PSOD or MSOD, $\hat{\tau}_{kk'} - \text{SATE}_{kk'} \to 0$ in probability.

### 3.3. Linear rate of convergence for parametric designs

In this section, we study the rate of convergence of the estimator under the PSOD or MSOD by studying the convergence of $\mathbb{E}[M_{\text{opt}}^2]$, which bounds it (see corollary 2). In particular, we argue that $\mathbb{E}[M_{\text{opt}}^2] = 2^{-\Omega(n)}$ when $m = 2$ and $\mathcal{F}$ finite dimensional. Empirically, $m \geqslant 3$ has similar convergence.

Let $\phi_1, \ldots, \phi_r$ be a basis for the finite dimensional $\mathcal{F}$ and $\Phi_{ij} = \phi_j(X_i)$. Because all norms in finite dimensions are equivalent, i.e. $c\|\cdot\|' \leqslant \|\cdot\| \leqslant C\|\cdot\|'$ (Hunter and Nachtergaele (2001), theorem. 5.36), it follows that any rate of convergence that applies when $\mathcal{F}$ is endowed with the 2-norm ($\|\beta_1\phi_1 + \ldots + \beta_r\phi_r\| = \|\beta\|_2$) also applies when $\mathcal{F}$ has any given norm. Moreover, note that, since $M_{\text{M-opt}}^2 \leqslant M_{\text{P-opt}}^2$, any rate for $\mathbb{E}[M_{\text{P-opt}}^2]$ applies also to $\mathbb{E}[M_{\text{M-opt}}^2]$. So, we restrict our attention to $\mathbb{E}[M_{\text{P-opt}}^2]$ for $\mathcal{F}$ under the 2-norm since any rate that applies to it will also apply generally.

Our argument is heuristic (not a precise proof) approximating the configurations $W$ with energies $M_P^2(W)$ as a spin glass following the random-energy model where energies are assumed independent. This approximation is commonly used to study the distributions of the optima of combinatorial optimization problems with random inputs and has been found to be valid asymptotically for similar partition problems (Mertens, 2001; Borgs *et al.*, 2009a, b).

Let $\Sigma_{ij} = \text{cov}\{\phi_i(X_1), \phi_j(X_1)\}$ and $\lambda_1, \ldots, \lambda_{r'} > 0$ be its positive eigenvalues where $r' = \text{rank}(\Sigma)$. The distribution of $M_P^2(W)$ is the same for any one fixed $W$. Fix $W_i = (i \bmod 2) + 1$ (i.e. $u_i = (-1)^{i+1}$). By the multivariate central limit theorem and continuous transformation, we have

$$\frac{2}{n}\Phi^\mathsf{T} u = \left( \frac{2}{n} \sum_{i=1}^{n/2} \{\phi_j(X_{2i-1}) - \phi_j(X_{2i})\} \right)_{j=1}^{r} \xrightarrow{\text{d}} \mathcal{N}(0, 2\Sigma),$$

$$M_P^2(W) = \sup_{\|\beta\|_2 \leqslant 1} \left\{ \frac{2}{n} \sum_{i=1}^{n} \sum_{j=1}^{r} u_i \beta_j \phi_j(X_i) \right\}^2 = \left\| \frac{2}{n} \Phi^{\mathsf{T}} u \right\|_2^2 \xrightarrow{\mathrm{d}} \sum_{i=1}^{r'} 2\lambda_i \chi_1^2,$$

the weighted sum of independent $\chi^2$ random variables with 1 degree of freedom. Denote the corresponding cumulative distribution function by $H$ and probability density function by $h$, which were given in series representation by Kotz *et al.* (1967). In following with the random-energy model approximation we assume independent energies so that $M_{P\text{-opt}}^2$ is distributed as the smallest order statistic among $\binom{n}{n/2}$-many independent draws from $H$. By Ahsanullah *et al.* (2013), theorem 11.3, and $\lim_{t \to 0^+} th(t)/H(t) = r'/2$, we have that $\mathbb{P}(M_{P\text{-opt}}^2/\beta_n \leqslant t) \to 1 - \exp(-t^{r'/2})$ for $\beta_n$ satisfying $\binom{n}{n/2} H(\beta_n) \to 1$. By Kotz *et al.* (1967), equation (40),

$$\beta_n = 4 \left\{ \frac{\Gamma(r'/2+1)}{\binom{n}{n/2}} \right\}^{2/r'} \prod_{i=1}^{r'} \lambda_i^{1/r'}.$$

Thus, $\mathbb{E}[M_{P\text{-opt}}^2] \approx \beta_n \Gamma(2/r'+1)$. By Stirling's formula,

$$\mathbb{E}[M_{M\text{-opt}}^2] \leqslant \mathbb{E}[M_{P\text{-opt}}^2] = O(2^{-2n/r'} n^{1/r'}) = 2^{-\Omega(n)}.$$

We verify this empirically for a range of cases; Fig. 4. We consider $m = 2, 3$, $X_i \sim \mathcal{N}(0, I_d)$ and $\phi_\theta(x) = s^{1-\Sigma_i \theta_i} \prod_{i=1}^{d} x_i^{\theta_i}$, $d = 1, 2, 3$, $r = \binom{d+s}{s}$ (all monomials up to degree $s$) for $s = 1, 2, 3$, and $q$-norms $1, 2$ and $\infty$. All exhibit linear convergence (note the log-scale) and are almost identical over $q$.

## 4. Algorithms for optimal design

In this section, we address how to compute optimal designs. For complete randomization, blocking, and pairwise-matched allocation with $m = 2$, how to do so is already clear. Here we address the other designs that arose from our framework and in particular kernel allocation. We solve for the PSOD by using mixed integer programming. Solving for the MSOD proves to be too difficult in practice so we provide heuristics based on semidefinite programming.

### 4.1. Optimizing pure strategies

The PSOD optimization problem can be written as

$$
\begin{aligned}
\sqrt{(} \min_{W \in \mathcal{W}} M_P^2(W)) = \min_{\lambda \in \mathbb{R}, w \in \{0,1\}^n} \ & \lambda \\
\text{subject to } \lambda \geqslant \max_{\|f\| \leqslant 1} \ & \frac{1}{p} \sum_{i=1}^{n} (w_{ik} - w_{ik'}) f(X_i) \quad \forall k < k' \\
\sum_{k=1}^{m} w_{ik} = 1 \ \forall i = 1, \dots, n, \quad & \sum_{i=1}^{n} w_{ik} = p \ \forall k = 1, \dots, m.
\end{aligned}
\tag{4.1}
$$

Next we show how to write constraints (4.1) to yield a mixed integer linear programme (MILP), mixed integer quadratic programme or a mixed-integer second-order cone programme (MIS-OCP). Since solver software may not return any one optimal solution at random, we randomly permute the optimal partition to ensure that assumption 1, part (b), holds. In our experiments, we use Gurobi version 5.6 (Gurobi Optimization, 2015).

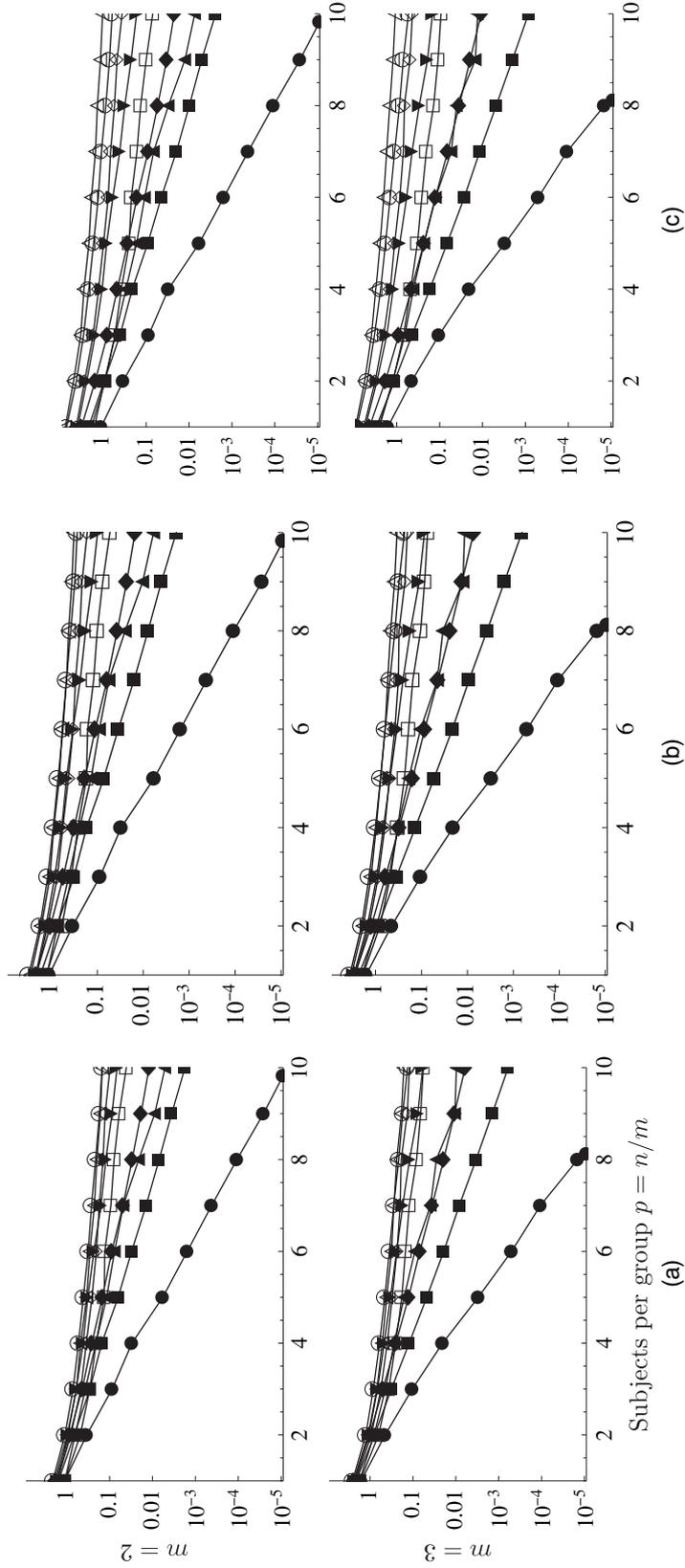

**Fig. 4.**  Convergence of $\mathbb{E}[M_{s,pl}^2]$ for $\binom{d+s}{s}$-dimensional $\mathcal{F}$ ($\bullet$, $d=1$, $s=1$; $\blacksquare$, $d=2$, $s=1$; $\blacklozenge$, $d=3$, $s=1$; $\blacktriangle$, $d=1$, $s=2$; $\blacktriangledown$, $d=2$, $s=2$; $\bigcirc$, $d=3$, $s=2$; $\square$, $d=1$, $s=3$; $\Diamond$, $d=2$, $s=3$; $\triangle$, $d=3$, $s=3$): (a) $q=1$; (b) $q=2$; (c) $q=\infty$

### 4.1.1. Finite dimensional q-space

For the set-up as in Section 2.3.4 and for $1/q + 1/q^* = 1$, we have

$$\max_{\|f\| \leqslant 1} \frac{1}{p} \sum_{i=1}^{n} (w_{ik} - w_{ik'}) f(X_i) = \left\| \left\{ \frac{n}{2} \sum_{i=1}^{n} (w_{ik} - w_{ik'}) \phi_1(X_i), \ldots, \frac{n}{2} \sum_{i=1}^{n} (w_{ik} - w_{ik'}) \phi_r(X_i) \right\} \right\|_{q^*}.$$

It follows that, for $q = 1, \infty$, constraint (4.1) is linear and the PSOD problem is an MILP. For $q = 2$, the problem for $m = 2$ is a mixed integer quadratic programme and for $m \geqslant 3$ it is an MISOCP (the difference being whether the quadratic term is in the objective or constraints). Finally, rational $q$ also leads to an MISOCP via the results of Lobo *et al.* (1998).

Since solving MILPs is generally faster than solving MISOCPs, it may be preferable to use the 1- or $\infty$-norms purely for the lower computational burden, especially for very large $n$. Either of these norms differs from the 2-norm by at most a factor of $\sqrt{r}$, and, as we saw in Fig. 4, the value and rate of convergence of $\mathbb{E}[M_{\text{P-opt}}^2]$ under these various $q$-norms is nearly identical.

### 4.1.2. Lipschitz functions

Next we consider the norm $\|f\| = \|f\|_{\text{lip}}$. When $m = 2$, theorem 4 shows that the PSOD is pairwise-matched allocation. The corresponding optimization problem is weighted non-bipartite matching, which is solved in polynomial time by using Edmond's algorithm (Edmonds, 1965). For $m \geqslant 3$, we let $D_{ij} = \delta(X_i, X_j)$ and use linear optimization duality (Bertsimas and Tsitsiklis, 1997) to obtain

$$\lambda \geqslant \max_{\|f\| \leqslant 1} \frac{1}{p} \sum_{i=1}^{n} (w_{ik} - w_{ik'}) f(X_i) = \max_{ve^{\mathsf{T}} - ev^{\mathsf{T}} \leqslant D} \frac{1}{p} \sum_{i=1}^{n} (w_{ik} - w_{ik'}) v_i$$

$$\Leftrightarrow \exists\, S \in \mathbb{R}_+^{n \times n} \quad \text{such that } \lambda \geqslant \text{tr}(DS)/p, \; \sum_{j=1}^{n} (S_{ij} - S_{ji}) = w_{ik} - w_{ik'} \quad \forall\, i = 1, \ldots, n,$$

yielding an MILP for the design problem.

### 4.1.3. Reproducing kernel Hilbert space

For RKHS norms, theorem 6 gives

$$\left\{ \max_{\|f\| \leqslant 1} \frac{1}{p} \sum_{i=1}^{n} (w_{ik} - w_{ik'}) f(X_i) \right\}^2 = \frac{1}{p} \sum_{i,j=1}^{n} (w_{ik} - w_{ik'}) K_{ij} (w_{jk} - w_{jk'}).$$

Therefore, for $m \geqslant 3$ the PSOD problem is an MISOCP and, for $m = 2$, we obtain the binary mixed integer quadratic programme

$$M^2(W) = \frac{4}{n^2} \min_{u \in \mathcal{U}} u^{\mathsf{T}} K u. \tag{4.2}$$

### 4.2. Optimizing mixed strategies

For the case of mixed strategies we consider only the case of $m = 2$ and $\mathcal{F}$ being an RKHS. As per theorems 2 and 6, the corresponding optimization problem is

$$\frac{4}{n^2} \min_{Q \in \mathcal{Q}} \lambda_{\max}(\sqrt{K} Q \sqrt{K}). \tag{4.3}$$

Although expression (4.3) has a convex objective and convex (polytope) feasible region, we have already observed in Section 2.4 that the problem is NP hard. But what makes problem (4.3) more difficult than problem (4.2) in practice is that it is not amenable to the branch-and-bound techniques that are employed by integer optimization software because its optimum generally does not occur at a corner point of the polytope, as we observed in Section 2.4. Therefore, we propose only heuristic solutions to the problem based on semidefinite programming, optimization over the cone $S_+^n$ of $n \times n$ positive semidefinite matrices (see Boyd and Vandenberghe (2004)). In our experiments, we use MOSEK version 7 (MOSEK, 2015) to solve semidefinite programs. The first heuristic is based on a semidefinite outer approximation of $\mathcal{Q}$ and the second on a simplicial inner approximation.

*Algorithm 1.* Let $\hat{Q}$ be a solution to the semidefinite program

$$\min_{\lambda \in \mathbb{R}, Q \in S_+^n} \lambda \qquad \text{subject to } \lambda I - \sqrt{K} Q \sqrt{K} \in S_+^n, \quad \text{diag}(Q) = e, \quad Qe = 0.$$

Let $\hat{\sigma}$ be the distribution of $u_i = \text{sgn}\{v_i - \text{median}(v)\}$ where $v \sim \mathcal{N}(0, \hat{Q})$.

*Algorithm 2.* Given $u_1, \ldots, u_T \in \mathcal{U}$, let $\hat{\theta}$ be a solution to the semidefinite program

$$\min_{\lambda \in \mathbb{R}, \theta \in \mathbb{R}^T} \lambda \qquad \text{subject to } \lambda I - \sum_{t=1}^{T} \theta_t \sqrt{K} u_t u_t^{\text{T}} \sqrt{K} \in S_+^n, \quad \theta \geqslant 0, \quad \sum_{t=1}^{T} \theta_t = 1.$$

Let $\hat{\sigma}$ be the distribution of $u = \pm u'$ where the sign is chosen equiprobably and where $u'$ is drawn randomly from $u_1, \ldots, u_T$ according to weights $\hat{\theta}_1, \ldots, \hat{\theta}_T$.

The inputs to algorithm 2 can be generated by running algorithm 1 and drawing some $u_t$ from the solution or by taking the top $T$ solutions to the PSOD problem as below. We use this in our experiments.

*Algorithm 3.* Set $\mathcal{U}_1 = \mathcal{U} \cap \{u_1 = 1\}$. Solve $u_t \in \arg\min_{u \in \mathcal{U}_t} u^{\text{T}} K u$, and set

$$\mathcal{U}_{t+1} = \mathcal{U}_t \cap \{u_t^{\text{T}} u \leqslant n - 4\}$$

for $t = 1, \ldots, T$. Run algorithm 2 using $u_1, \ldots, u_T$.

## 5. Algorithms for inference

Since balance can significantly reduce estimation variance, we would expect that it can also increase statistical power. In this section, we consider $m = 2$ and the testing of the sharp null hypothesis

$$H_0 : (\text{TE}_i = 0 \ \forall \ i = 1, \ldots, n).$$

Under $H_0$ all $n$ outcomes are exchangeable regardless of treatment given ($Y_{i1} = Y_{i2}$). We can therefore simulate what would happen under another assignment and compare. This is the idea behind Fisher's randomization test, where new simulated assignments are drawn from the same design as used at the outset of the experiment. Fisher's randomization test can be applied to *any a priori* balancing design and it will always yield an exact $p$-value and provably has type I error rate no greater than the designed significance. This test is standard and we defer discussion on its effective use for MSODs and the outputs of algorithms 1 and 2 to the on-line supplemental section D.

A special case occurs for *a priori* balancing designs that only randomize over treatment permutations of a single partition of subjects, such as rerandomization with infinitesimal acceptance probability or kernel allocation PSODs. Although Fisher's randomization test is valid (correct significance), it will also have zero power against any alternative (it always yields $p$-value 1 because all permutations will yield a statistic of the same magnitude). To address this deficiency, we develop an alternative test based on the bootstrap, which we believe is asymptotically valid (Efron and Tibshirani, 1993).

*Algorithm 4.* For a confidence level $0 < 1 - \alpha < 1$, perform the following steps.

*Step 1*: draw $W^0$ from the PSOD for the baseline covariates $X_1, \ldots, X_n$, assign subjects, apply treatments, measure outcomes $Y_i = Y_{iW_i^0}$ and compute $\hat{\tau}$.
*Step 2*: for $t = 1, \ldots, T$

    (a) sample $i_j^t \sim \mathrm{unif}(1, \ldots, n)$ independently for $j = 1, \ldots, n$,
    (b) draw $W^t$ from the PSOD for the baseline covariates $X_{i_1^t}, \ldots, X_{i_n^t}$ and
    (c) compute $\tilde{\tau}^t = (1/p)\Sigma_{j:W_j^t=1} Y_{i_j^t} - (1/p)\Sigma_{j:W_j^t=2} Y_{i_j^t}$.

*Step 3*: the $p$-value of $H_0$ is $p = (1 + |\{t : |\tilde{\tau}^t| \geqslant |\hat{\tau}|\}|)/(1 + T)$. If $p \leqslant \alpha$, then reject $H_0$.

Algorithm 4 can also be used as a test for any *a priori* balancing design including the MSOD, letting the new design be computed for the bootstrap sample in step 2(b). However, the additional randomization of the MSOD (and complete randomization, etc.) means that we do not have to resort to a bootstrap approximation as we can use the standard randomization test, which also eliminates the computational burden of having to recompute the design. See the supplemental section D for a discussion.

Using the duality of hypothesis tests and confidence intervals, algorithm 4 can also be used to construct confidence intervals around the point estimate $\hat{\tau}$ (see Good (2005), section 3.3).

## 5.1. Example 4

Consider the set-up as in example 2 with $d = 2$, quadratic $\hat{f}$ and $\epsilon_{i1} = \epsilon_{i2} \sim \mathcal{N}(0, \sigma_{\mathrm{residual}}^2)$. For various values of $\tau$ and $\sigma_{\mathrm{residual}}$, we test $H_0$ at significance $\alpha = 0.05$. We replace each kernel allocation MSOD in example 2 with the corresponding kernel allocation PSOD and use algorithm 4. For all other designs we use the standard randomization test (see algorithm D.2 in the on-line supplement). Fig. 5 shows the probability of rejecting $H_0$.

When $\tau$ is non-zero, the null hypothesis is false and the plot shows *power*, or the complement of the type II error rate. Quadratic and exponential kernel allocation detect the positive effect almost *immediately* and Gaussian kernel allocation soon after. The classical designs lag far behind in detection power. Linear kernel allocation parametrically misspecifies the conditional expectation function in this case but, nonetheless, does not do much worse than the best of the classical designs. As the magnitude of residual noise $\sigma_{\mathrm{residual}}$ increases, the estimation variance under any one design increases (by the same amount) as per theorem 7 and the relative gains in power due to balancing are slowly diminished.

When $\tau$ is 0, the null hypothesis is true and the plot shows the type I error rate. Since all type I error rates are no greater than 0.05, all the tests are valid at the designed significance $\alpha = 0.05$. However, the tests that are carried out by algorithm 4 for the kernel allocation PSODs are particularly conservative. As the imbalance disappears, these tests have a much lower type I error than the significance $\alpha = 0.05$. We emphasize that this is a *desirable* aberration—this means that these tests have *both* incredibly high power and extra low type I error, i.e. they strictly dominate the other tests. This conservatism does suggest, however, that we may be able to achieve even higher power while maintaining significance of $\alpha = 0.05$ by artificially deflating

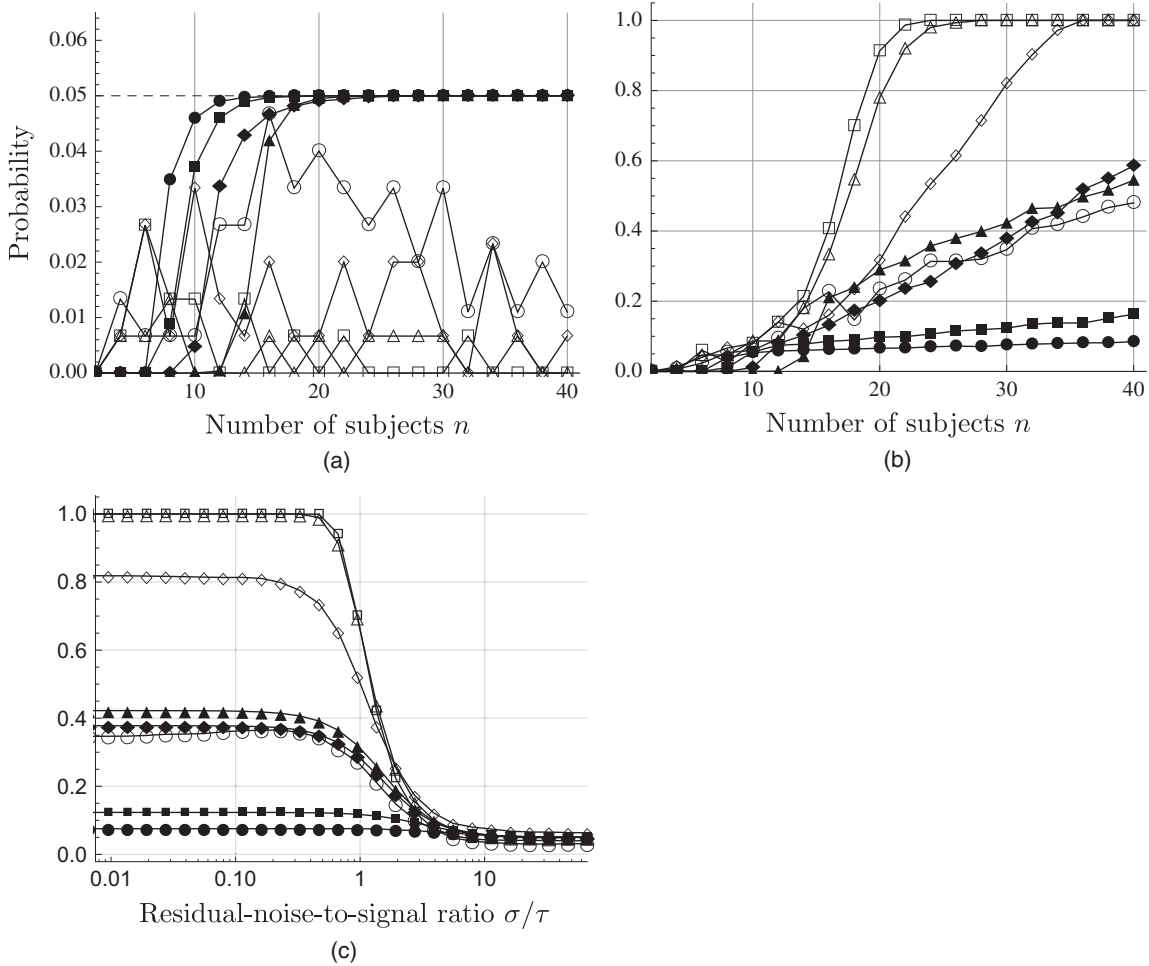

**Fig. 5.** Probability of rejecting $H_0$ at $\alpha = 5\%$ under various designs as in example 4 and no effect ($\tau = 0$), a positive effect ($\tau = 0.15$) or a positive effect with varying residual noise ($\bullet$, optimal $L^{\infty}$ (complete randomization); $\bigcirc$, optimal linear; $\blacksquare$, optimal $L^{\infty}$ coarsened (blocking); $\square$, optimal quadratic; $\blacklozenge$, optimal Lipschitz (pairwise matching); $\lozenge$, optimal Gaussian; $\blacktriangle$, rerandomization; $\triangle$, optimal exponential): (a) $\tau = 0$, $\sigma = 0$, $n$ varies; (b) $\tau = 0.15$, $\sigma = 0$, $n$ varies; (c) $\tau = 0.15$, $\sigma$ varies, $n = 30$

the critical $p$-value in step 3 of algorithm 4. (This conservatism artefact does not influence the application of the standard randomization test to kernel allocation MSODs; here we are using kernel allocation PSODs. A possible alternative for testing under the PSOD may be to take the top 20 solutions to the PSOD optimization problem, to randomize among these and to use Fisher's randomization test. This increases type I error to be exactly 0.05 but also reduces power by giving equal randomization weight to suboptimal solutions: the power of this test for quadratic kernel allocation is reduced from about 100% to about 90% at $n = 24$, i.e. this is *strictly* dominated by our bootstrap approach.)

## 6.  General recommendations for practice

In this paper, we developed a general framework for optimal *a priori* balance that encompasses many existing designs and also gives rise to new, powerful kernel-based optimal designs for controlled experiments. The totality of the results in this paper—theoretical and experimental— suggests that it is a very rare case where kernel allocation should not be favoured for use in practice.

If the baseline covariates are irrelevant to begin with, then any design does equally well and there is no way to improve or do worse (theorem 9). If the baseline covariates are relevant and the true conditional expectation function is well specified or well approximable by $\mathcal{F}$, then balancing the covariates optimally can lead to nearly optimal reduction in variance for even moderate $n$ (theorem 8, theorem 10 and corollary 3, section 3.3). In particular, for universal kernel allocation, we can arbitrarily closely approximate any function, leading to model-free consistency (theorems 11 and 13). Usually, non-linear polynomials are sufficient for good approximation in practice. Dealing with clinical data, we saw in example 3 that quadratics offered a successful balance between generality of structure and efficiency of balancing all possible quadratic functions simultaneously. In example 2, we saw that even rather extreme violations of quadraticity did not significantly impede performance. If, however, the true conditional expectation function is completely perverse relative to assumed structure then we are in the setting of theorem 1 and may end up in a situation like that in example 1, where outcomes perversely depended on the parity of the logarithm of the covariates and any balancing (even by classical designs) only hurts precision. One might argue that such examples are of limited relevance in practice.

Although polynomial kernel allocation can offer good performance in practice, one drawback and potential criticism is that it does not guarantee consistency if the true conditional expectation function is not polynomial (of bounded degree). This is in contrast with complete randomization and blocking, which are always consistent under mild conditions (existence of moments). However, universal kernel allocation is also always consistent (see theorem 13) and, in every single example that we considered, universal kernel allocation outperformed each of complete randomization, blocking, pairwise-matched allocation and rerandomization. Therefore, if consistency is a concern beyond small sample efficiency, then universal kernel allocation should be used instead of polynomial kernel allocation.

When using any kernel that is not translation invariant (like the Gaussian kernel is), it is advisable to centre the data by subtracting the mean. When the different covariates are not comparable (e.g. measured in different units), it is advisable to rescale the centred data either by dividing by standard deviations or by premultiplying by the inverse square root of the sample covariance as in example 3, which ensures that kernel allocation with a kernel that depends only on inner products or distances (this includes all kernels introduced here) will be invariant to affine transforms of the data. It should be noted that our theoretical results deal with independent subjects so this practical stopgap does not fit exactly into this framework (the covariate vectors are no longer independent after recentring). However, this issue becomes moot for $n$ sufficiently larger than the dimension $d$ of the covariate data.

## 7.   Concluding remarks

Designs that provide balance in controlled experiments before treatment and before randomization provide one answer to the criticism that complete randomization may lead to assignments that the experimenter knows will lead to misleading conclusions. In this paper we unified these designs under the umbrella of *a priori* balance. We argued that structural information on the dependence of outcomes on baseline covariates was the key to any *a priori* balance beyond complete randomization and developed a framework of optimal designs based on structure expressed on the conditional expectation function. We have shown how existing *a priori* balancing designs, including blocking, pairwise-matched allocation and other designs, are optimal for certain structures and how existing imbalance metrics, such as Mahalanobis distance between group means, arise from other choices of structure. That this theoretical framework fitted so well into existing practice led us to endeavour to discover what other designs may arise from it. We considered a wide range of designs that follow from structure expressed by using kernels, encompassing both parametric and non-parametric methods. We argued and showed numerically that parametric models (when correctly specified) coupled with optimization lead to estimation variance that converges very fast to the best theoretically possible. The implication in practice, as seen in experiments with clinical data, is that such designs can successfully leverage additional, even if marginal, prognostic content in additional dimensions of baseline data, whereas many classic designs struggle with even a few additional dimensions. The new connections that were made in this paper suggest many avenues for future research, including studying designs arising from new choices of norms, convergence rates of minimal imbalance for infinite dimensional spaces, data-driven choices of structure by, for example, regression analysis on initial data, optimal *a priori* balance in sequential design where allocations must be made on the fly, further study of the bootstrap testing procedure and its validity, and specialized algorithms for solving the design optimization problems.

## Acknowledgements

I thank the reviewers and Associate Editor for their extremely helpful suggestions and very thorough review of the paper.

*Supporting information*
Additional 'supporting information' may be found in the on-line version of this article:
  'Supplement to: Optimal a priori balance in the design of controlled experiments',

**Supplement to:**

**Optimal A Priori Balance in the Design of Controlled Experiments**

Nathan Kallus

*Cornell University, New York, NY, USA.*

## A. A priori balance in estimating treatment effect on compliers

In many experimental endeavors involving human subjects the researcher does not fully control the treatment actually administered. Consider two treatments, "treatment" ($k = 1$) and "control" ($k = 2$). Situations where a subject receives a treatment different from their assignment include refusal of surgery, ethical codes that allow subjects assigned to control to demand treatment, or the leakage of information to some control subjects in a teaching intervention. This issue is termed non-compliance. In such situations, $W$ represents initial assignment intent and our estimator $\hat{\tau}$ estimates the effect of the *intent* to treat (ITT). Often a researcher is interested in the compliers' average treatment effect in the sample (CSATE) or population (CPATE), disregarding all non-compliers. Subjects that always demand treatment are known as always-takers, those that always refuse treatment as never-takers, and those that always choose the opposite of their assignment as defiers (this is exhaustive if subjects comply based only on their own assignment). Denote by $\pi_c$ and $\Pi_c$ the unknown fraction of compliers in the sample and population, respectively. In the absence of defiers we can observe the identity of never-takers in the treatment group and of always-takers in the control group. We can estimate the fraction of compliers as the complement of those:

$$\hat{\pi}_c = 1 - \frac{2}{n} \sum_{i:W_i=1} \mathrm{NT}_i - \frac{2}{n} \sum_{i:W_i=2} \mathrm{AT}_i$$

where $\mathrm{NT}_i = 1$ if $i$ is a never-taker and $\mathrm{AT}_i = 1$ if $i$ is an always-taker (both 0 for compliers). Under an assignment that blinds the identity of treatment, such as complete randomization, $\hat{\pi}_c$ is conditionally (for $\pi_c$) and marginally (for $\Pi_c$) unbiased if there are no defiers. Without defiers,

$$\mathrm{CSATE} = \mathrm{SATE} / \pi_c \qquad \mathrm{CPATE} = \mathrm{PATE} / \Pi_c$$

since the individual ITT effect for an always- or never-taker is identically 0. The standard approach in completely randomized trials is to estimate the compliers' average treatment effect by the Wald ratio estimator $\hat{\tau}_c = \hat{\tau}/\hat{\pi}_c$ (Imbens and Rubin, 1997; Little and Yau, 1998). Such an estimator, equivalent to the two-stage least squares estimate (without covariates), need not be unbiased.

We can do even better if we use a priori balance to improve the precision of the compliance fraction estimator. The difference between the sample compliance fraction and our estimator of it is

$$\hat{\pi}_c - \pi_c = \tfrac{2}{N} \sum_{i:W_i=1} (\text{AT}_i - \text{NT}_i) - \tfrac{2}{n} \sum_{i:W_i=2} (\text{AT}_i - \text{NT}_i) = \tfrac{2}{n} \sum_{i=1}^n u_i C_i$$

$$\text{where} \quad C_i = \left\{ \begin{array}{ll} 1 & i \text{ is always-taker} \\ 0 & i \text{ is complier} \\ -1 & i \text{ is never-taker} \end{array} \right. \quad \text{is } i\text{'s compliance status.}$$

Therefore, matching the means of $f_c(x) = \mathbb{E}[C_i \mid X_i = x]$ will eliminate variance in estimating the compliance fraction and get us closer to the true CSATE and CPATE. Moreover, if the two unbiased estimators, $\hat{\tau}$ and $\hat{\pi}_c$, are both more precise, their ratio $\hat{\tau}_c$ is both more precise and less biased. To achieve this through our framework we need only incorporate our belief $\mathcal{F}_c$ about $f_c$ into the larger $\mathcal{F}$ and proceed as before. (See also supplemental Sec. B for a discussion about combining spaces.)

## B.  Generalizations of $\mathcal{F}$

In this supplemental section we consider more general forms of the space $\mathcal{F}$. For the most part, the theorems presented in the main text will still apply. We deferred this discussion to this supplement to avoid overly cumbersome notation in the main text.

First, we consider the restriction to cones in $\mathcal{F}$. A cone is a set $C \subset \mathcal{F}$ such that $f \in C \implies cf \in C \ \forall c > 0$. We may then further restrict to $f \in C$, $||f|| \leq 1$ in the definitions of $M_{\text{P}}^2(W)$ and $M_{\text{M}}^2(\sigma)$. By symmetry, this is the same as restricting to $C \cup (-C)$. Since it is still the case that $||cf|| = c||f||$, Thms. 8 and 12 still apply. One example of a cone is the cone of monotone functions (either nondecreasing or nonincreasing). In a single dimension and for two treatments, this will result in a PSOD that sorts the data and assigns subjects in an alternating fashion. This is also an optimal assignment for pairwise-matched allocation in one dimension. More generally and in higher dimensions, we can consider a directed acyclic graph (DAG) on the nodes $V = \{1, \ldots, n\}$ with edge set $E \subset V^2$ and its associated topological cone $C = \{f : f(X_i) \leq f(X_j) \ \forall (i,j) \in E\}$. Other cones include nonnegative/positive functions and $\pm$-sum-of-squares polynomials.

Second, we consider re-centering the norms. We might have a nominal regression function $g$ that we believe is approximately right, perhaps due to a prior regression analysis or based on models from the literature. In this case, it would make sense to solve the minimax problem against perturbations around this $g$. Given a norm $||\cdot||'$ on $\mathcal{F}$ we can formally define the magnitude

$$||f|| = \max \left\{ \min \left\{ ||f - g||', ||f + g||' \right\}, 1 \right\}. \tag{B.1}$$

We consider both $g$ and $-g$ because it has no effect on the imbalance metrics due to symmetry of the objective while it can only reduce magnitudes. Using this alternate definition of $||\cdot||$ in (B.1), Thm. 8 still applies and Thm. 12 applies if its conditions apply to the Banach space $\mathcal{F}$ with its usual norm and $\mathbb{E} |g(X_1)| < \infty$. In the Bayesian interpretation discussed in Sec. 2.4.1, this is equivalent to making the prior mean of $f(x)$ be $g(x)$.

Third, we consider combining multiple spaces $\mathcal{F}_1, \ldots, \mathcal{F}_b$. There are two ways. The first way is to combine these via an algebraic sum. The space $\mathcal{F} = \mathcal{F}_1 + \cdots + \mathcal{F}_b = \{\phi_1 + \cdots + \phi_b : \phi_j \in \mathcal{F}_j \ \forall j\}$ endowed with the norm $||f|| = \min_{\phi_j \in \mathcal{F}_j : f = \phi_1 + \cdots + \phi_b} \max_{j=1,\ldots,b} ||\phi_j||_{\mathcal{F}_j}$ is a Banach space and as such a valid choice. In particular, the algebraic sum $\mathcal{F}$ can be identified with the quotient of the direct sum $\mathcal{F}' = \mathcal{F}_1 \oplus \cdots \oplus \mathcal{F}_b$ by its subspace $\{(\phi_1, \ldots, \phi_b) \in \mathcal{F}' : \phi_1 + \cdots + \phi_b = 0\}$. We can decompose the pure-strategy imbalance metric corresponding to this new choice as follows:

$$M_{\mathrm{P}}^2(W) = \max_{k \neq k'} \left( \sum_{j=1}^{b} \sup_{||\phi_j||_{\mathcal{F}_j} \leq 1} B_{kk'}(W, \phi_j) \right)^2.$$

Thms. 8 and 12 still apply and the conditions of Thm. 12 hold for $\mathcal{F}$ if they hold for each $\mathcal{F}_j$.

The second way is to combine these formally via a union. Consider the space $\mathcal{F} = \mathcal{F}_1 \cup \cdots \cup \mathcal{F}_b = \{f : f \in \mathcal{F}_j \text{ for some } j\}$. This is not a vector space but we can formally define the magnitude $||f|| = \min_{j=1,\ldots,b} ||f||_{\mathcal{F}_j}$. We can then decompose the pure-strategy imbalance metric corresponding to this new choice as follows:

$$M_{\mathrm{P}}^2(W) = \max_{k \neq k'} \max_{j=1,\ldots,b} \sup_{||\phi_j||_{\mathcal{F}_j} \leq 1} B_{kk'}^2(W, \phi_j).$$

Thm. 8 still applies and Thm. 12 applies if its conditions hold for each Banach space $\mathcal{F}_j$.

We can even take several spaces $\mathcal{F}_1, \ldots, \mathcal{F}_b$, re-center each norm with its own $g_j$ as in (B.1), and then combine them in either of the two ways, defining the combined magnitudes strictly formally. In this way, we can have multiple centers to represent various beliefs about the same or different regression functions $f_k$. Thm. 8 still applies and Thm. 12 applies if its conditions hold for each $\mathcal{F}_j$ and $\mathbb{E}\,|g_j(X_1)| < \infty$ for for each $j$.

## C. Guarantees on general moment mismatch

In Sec. 3.1, we obtained strong guarantees on the estimation variance of optimal a priori balancing designs by bounding the mismatch in the generalized moment of $f_k$ between the experimental groups. We can similarly bound the mismatch in any other moment by plugging in a different function.

For example, we may wish to derive a bound on the mismatch in first moments (where $\mathcal{X} \subset \mathbb{R}^d$). For any $v \in \mathbb{R}^d$, define the function $\phi_v(x) = v^T x$. Then, for any positive semi-definite matrix $Q$, let

$$\mathbb{E}M_Q^2 := \mathbb{E}[(\tfrac{2}{n}\sum_{i=1}^{n} u_i X_i)^T Q \left(\tfrac{2}{n}\sum_{i=1}^{n} u_i X_i\right)] = \mathbb{E}[\sup_{||v|| \leq 1} \left(v^T \sqrt{Q}\left(\tfrac{2}{n}\sum_{i=1}^{n} u_i X_i\right)\right)^2]$$
$$= \mathbb{E}[\sup_{||v|| \leq 1} B^2(W, \phi_{\sqrt{Q}v})] \leq \mathbb{E}M_{\mathrm{opt}}^2 \sup_{||v|| \leq 1} ||\phi_{\sqrt{Q}v}||^2.$$

So to produce a bound, we need only calculate $\sup_{||v|| \leq 1} ||\phi_{\sqrt{Q}v}||^2$. If, for example, we use the quadratic kernel then we can write $\phi_v(x) = v^T x = \frac{1}{2}(1 + v^T x)^2 - \frac{1}{2}(1 - v^T x)^2$ and deduce that

$$||\phi_v||^2 = \begin{pmatrix} 1/2 \\ 1/2 \end{pmatrix}^T \begin{pmatrix} (1 + v^T v/2)^2 & (1 - v^T v/2)^2 \\ (1 - v^T v/2)^2 & (1 + v^T v/2)^2 \end{pmatrix} \begin{pmatrix} 1/2 \\ 1/2 \end{pmatrix} = ||v||_2^2.$$

Similarly, it can be shown that in any polynomial RKHS and in the exponential RKHS, we also have $||\phi_v||^2 = ||v||_2^2$. Therefore, in all of these cases, we have

$$\mathbb{E}M_Q^2 \leq \mathbb{E}M_{\mathrm{opt}}^2 \sup_{v^T v \leq 1} ||\sqrt{Q}v||^2 = ||Q||_2 \,\mathbb{E}M_{\mathrm{opt}}^2, \tag{C.1}$$

where $||Q||_2$ is the operator norm (largest singular value). If $Q$ is the identity matrix ($||Q||_2 = 1$), then $\mathbb{E}M_Q^2$ is the sum of all squared mismatches.

Estimation variance is the target quantity to reduce when estimation is unbiased, and in this paper we argued that the appropriate metric for imbalance is worst-case variance. However, it is a common practice to check the imbalance of a given design with respect to mismatch in means even if one does not assume a linear effect. These bounds provide a guarantee on this mismatch.

## D. Inference for mixed-strategy designs

As noted in Sec. 5, Alg. 4 can be used to answer inferential questions for mixed-strategy designs as well, but their additional randomization allows for the standard randomization and exact permutation tests to be used instead. The following is the standard permutation test when applied to a non-completely randomized design, including the MSOD.

*Algorithm* 1. Let $\sigma$ be given. For a confidence level $0 < 1 - \alpha < 1$:

1: Draw $W^0$ from $\sigma$, assign subjects, apply treatments, measure $Y_i = Y_{iW_i^0}$, and compute $\hat{\tau}$. Let $\mathcal{W}' = \{W \in \mathcal{W} : \sigma(W) > 0\}$.

2: For $W \in \mathcal{W}'$ compute $\tilde{\tau}^W = \frac{1}{p}\sum_{i:W_i=1} Y_i - \frac{1}{p}\sum_{i:W_i=2} Y_i$.

3: The $p$-value of $H_0$ is $p = \sum_{W \in \mathcal{W}'} \sigma(W)\mathbb{I}\left[\left|\tilde{\tau}^W\right| \geq |\hat{\tau}|\right]$. If $p \leq \alpha$ then reject $H_0$.

Note that for good power with the randomization we should have $\sigma(\{W, -W\}) \leq \alpha$ for every $W \in \mathcal{W}'$ (in particular, this is not true for the PSOD). To guarantee this for the output of Alg. 2, we can constrain the mixture components $\theta$ to be no greater than $\alpha$ in each component.

The above exact test requires that we have a full description of $\sigma$ and that we iterate over all feasible assignments. This works well for the output of Alg. 2 but can be prohibitive for the output of Alg. 1. The standard randomization test eschews these issues.

*Algorithm* 2. Let $\sigma$ be given. For a confidence level $0 < 1 - \alpha < 1$:

1: Draw $W^0$ from $\sigma$, assign subjects, apply treatments, measure $Y_i = Y_{iW_i^0}$, and compute $\hat{\tau}$.

2: For $t = 1, \ldots, T$ do:

  2.1: Draw $W^t$ from $\sigma$.

  2.2: Compute $\tilde{\tau}^t = \frac{1}{p}\sum_{i:W_i^t=1} Y_i - \frac{1}{p}\sum_{i:W_i^t=2} Y_i$.

3: The $p$-value of $H_0$ is $p = \left(1 + \left|\left\{t \ : \ \left|\tilde{\tau}^t\right| \geq |\hat{\tau}|\right\}\right|\right) / (1 + T)$. If $p(H_0) \leq \alpha$ then reject $H_0$.

### D.1. Aside: validity

A design is termed "valid" if the ensuing estimation or testing based upon it is valid (Moulton, 2004; Bailey, 1983; Bailey and Rowley, 1987). As such, as pointed out by Bailey and Rowley (1987), validity depends upon the purpose of the study and the intentions of the experimenter. A basic requirement, common to nearly all notions of validity, is that the estimate for treatment effect (or other contrasts of interest) be unbiased by requiring the blinding condition in Asmp. 1(b). In particular, this is equivalent to the notion of validity of Kupper et al. (1981), which discusses validity in matched designs, and to the first-order validity conditions of White and Welch (1981). In this sense, any design that satisfies Asmp. 1(b) is valid because it leads to valid estimation of treatment effect. But, if testing is intended, the validity of a design is the question of whether it leads to valid testing.

Many works dealing with validity of design either assume a testing procedure based on asymptotic distributions ($t$-tests and $F$-tests) via estimates of standard errors or are seeking valid estimation of standard error for its own right (Bailey, 1983; White and Welch, 1981; Youden, 1972; Grundy and Healy, 1950). As per Kempthorne (1966), Fisher's main justification for randomization was that it guarantees "a valid estimator of error" and "the validity of the test of significance." In these works, validity is predicated on valid estimation of standard error, which of course in turn depends on how one estimates error. For example, the second-order validity conditions of White and Welch (1981) require that every pair of subjects have the same probability of receiving any pair of treatments. These conditions are necessary for the pooled sample variance estimator to be valid. White and Welch (1981); Youden (1972); Bailey (1983); Grundy and Healy (1950) and others consider restrictions on complete randomization that preserve this validity by satisfying these or similar conditions. However, these conditions are also clearly violated by standard designs like randomized block designs or any pair-matched allocation. But, in such designs, other estimators of the error are in fact appropriate, enabling an analysis tailored for the design (Fisher, 1935; Snedecor and Cochran, 1989; Abadie and Imbens, 2008). (Note that, in light of Thm. 2, variance of effect estimation must be exactly the same among all designs that satisfy the validity conditions of White and Welch, 1981.)

In this study, we present allocations that violate the second-order conditions of White and Welch (1981). However, we make no attempt to understand the asymptotic distribution of our estimator nor its standard error and thus are not bound by these or similar requirements. Instead, since our designs satisfy Asmp. 1, we can achieve valid inference by using Fisher's randomization test (employing the same randomization design as used at the onset), as discussed above and in Sec. D.

A special situation arises in the case of a design that randomizes over too few partitions. In this case, a randomization test is still valid (i.e., has type I error rate no greater than the prescribed significance), but it may lack power. In the extreme case where the design fixes the partition and only randomize the treatment labels, Fisher's randomization test will always yield a $p$-value of 1. It is precisely for this reason that we introduce the bootstrap test in the previous section. It remains an open question under what conditions will the bootstrap test be asymptotically valid. Nonetheless, we observe that, empirically, we maintain validity (i.e., it is conservative in its type I errors, relative to the prescribed significance).

## E. Proofs

*Proof of Thm. 1.* Note that

$$\text{Var}(\hat{\tau}^{\text{CR}} \mid X, Y) = \frac{4}{n(n-1)} \left|\left|\overline{Y}\right|\right|_2^2 \quad \text{where} \quad \hat{Y}_i = \frac{Y_{i1} + Y_{i2}}{2}, \quad \hat{\mu} = \frac{1}{n} \sum_{i=1}^n \hat{Y}_i, \quad \text{and} \quad \overline{Y}_i = \hat{Y}_i - \hat{\mu}.$$

$$\hat{\tau} - \text{SATE} = \frac{2}{n} \sum_{i:W_i=1} \left(\frac{Y_{i1} + Y_{i2}}{2}\right) - \frac{2}{n} \sum_{i:W_i=2} \left(\frac{Y_{i1} + Y_{i2}}{2}\right) = \frac{2}{n} \sum_{i=1}^n u_i \hat{Y}_i.$$

By conditional unbiasedness, we have

$$\mathrm{Var}(\hat{\tau} \mid X, Y) = \mathbb{E}[(\hat{\tau} - \mathrm{SATE})^2 \mid X, Y] = \mathbb{E}\left[\left(\frac{2}{n}\sum_{i=1}^n u_i\hat{Y}_i\right)^2 \mid X, Y\right] = \sum_{W \in \mathcal{W}}\sigma(W)\left(\frac{2}{n}\sum_{i=1}^n u_i\hat{Y}_i\right)^2.$$

Thus, letting $V' = n(n-1)V/4$ our objective is to minimize

$$\max_{Y \in \mathbb{R}^{n \times 2} : \|\bar{Y}\|_2^2 = V'}\sum_{W \in \mathcal{W}}\sigma(W)\left(\frac{2}{n}\sum_{i=1}^n u_i\hat{Y}_i\right)^2, = \max_{Y \in \mathbb{R}^{n \times 2} : \|\bar{Y}\|_2^2 \le V'}\sum_{W \in \mathcal{W}}\sigma(W)\left(\frac{2}{n}\sum_{i=1}^n u_i\hat{Y}_i\right)^2,$$

because given a feasible $\hat{Y}$ in the above, replacing it with $\overline{Y} = \hat{Y} - \hat{\mu}$ is still feasible and has the same objective, and given a feasible $\overline{Y}$, replacing it with $\sqrt{V'} \times \overline{Y}/\left\|\overline{Y}\right\|_2$ is still feasible and only increases the objective.

Note that $\|\hat{Y}\|_2$ is a row-permutationally invariant seminorm on $Y \in \mathbb{R}^{n \times 2}$. Consider, more generally, any such seminorm and the objective

$$R(\sigma) = \max_{Y \in \mathbb{R}^{n \times 2} : \|Y\| \le \sqrt{V'}}\sum_{W \in \mathcal{W}}\sigma(W)\left(\frac{2}{n}\sum_{i=1}^n u_i\hat{Y}_i\right)^2.$$

Suppose $\sigma \in \Delta$ minimizes $R(\sigma)$. For $\pi \in S_n$ a permutation of $\{1, \ldots, n\}$, define $\sigma_\pi((W_1, \ldots, W_n)) = \sigma((W_{\pi(1)}, \ldots, W_{\pi(n)}))$. Then by the symmetry of $\|\cdot\|$, $R(\sigma) = R(\sigma_\pi)$, so $\sigma_\pi$ is also optimal. Next note that $R(\sigma)$ is a maximum over linear forms in $\sigma$ and is hence convex. Therefore, $\sigma^*(W) = \frac{1}{n!}\sum_{\pi \in S_n}\sigma_\pi(W)$ is also optimal. By construction we get $\sigma^*((W_1, \ldots, W_n)) = \sigma^*((W_{\pi(1)}, \ldots, W_{\pi(n)}))$ for any $\pi \in S_n$. Hence, $\sigma^*((W_1, \ldots, W_n))$ is constant for every $W \in \mathcal{W}$, and therefore $\sigma^*$ is complete randomization.   $\square$

*Proof of Thm. 2.* First note that by Asmp. 1(b), for any $i, j, k, k'$,

$$\sigma(\{W : W_i = W_j, W_i \in \{k, k'\}, W_j \in \{k, k'\}\}) = \frac{2}{m}\sigma(\{W : W_i = W_j\}),$$
$$\sigma(\{W : W_i \ne W_j, W_i \in \{k, k'\}, W_j \in \{k, k'\}\}) = \frac{2}{m}\frac{1}{m-1}\sigma(\{W : W_i \ne W_j\}).$$

Therefore, for any $k \ne k'$ and $i, j$,

$$\sum_{W \in \mathcal{W}}\sigma(W)(w_{ik} - w_{ik'})(w_{jk} - w_{jk'}) = \frac{2}{m}Q_{ij}(\sigma).$$

Therefore, by expanding the square and interchanging sums, for any $k \ne k'$ we have

$$\sum_{W \in \mathcal{W}}\sigma(W)B_{kk'}^2(W, f) = \frac{1}{p^2}\sum_{i,j=1}^n f(X_i)f(X_j)\sum_{W \in \mathcal{W}}\sigma(W)(w_{ik} - w_{ik'})(w_{jk} - w_{jk'})$$
$$= \frac{2}{pn}\sum_{i,j=1}^n Q_{ij}(\sigma)f(X_i)f(X_j). \qquad \square$$

*Proof of Thm. 3.* Let $\{x_1, \ldots, x_\ell\}$ be the set of values taken by the baseline covariates $X_1, \ldots, X_n$ ($\ell \le n$). Let an assignment $W$ be given. Let $\{i_1, i_1'\}, \ldots, \{i_q, i_q'\}$ denote a maximal perfect exact match across the two groups ($W_{i_j} = 1$, $W_{i_j'} = 2$, $X_{i_j} = X_{i_j'}$, and $q$ maximal) with $\{i_1'', \ldots, i_{q'}''\}$, $\{i_1''', \ldots, i_{q'}'''\}$ being the remaining unmatched subjects ($W_{i_{j'}''} = 1$, $W_{i_{j'}'''} = 2$, $X_{i_{j'}''} \ne X_{i_{j'}'''}$). For $i = 1, \ldots, \ell$, if there are more $x_i$'s in group 1 set $f'(x_i) = 1$ otherwise set $f'(x_i) = -1$. This $f'$ has $\|f'\|_\infty \le 1$ and hence

$$\max_{\|f\| \le 1}|B(W, f)| \ge |B(W, f')| = \frac{2}{n} \times q' \times 2 = 2 - \frac{4}{n}q.$$

At the same time, we have

$$\max_{\|f\| \le 1}|B(W, f)| = \max_{\|f\| \le 1}\left|\sum_{i=1}^n u_i f(X_i)\right|$$
$$\le \frac{2}{n}\sum_{j=1}^q \max_{\|f\| \le 1}\left|f(X_{i_j}) - f(X_{i_j'})\right| + \frac{2}{n}\sum_{j=1}^{q'}\max_{\|f\| \le 1}\left|f(X_{i_{j'}''}) - f(X_{i_{j'}'''})\right|$$
$$= 0 + \frac{2}{n} \times q' \times 2 = 2 - \frac{4}{n}q.$$

To summarize,

$$\sqrt{M_P^2(W)} = 2 - \frac{4}{n}\begin{pmatrix}\text{number of perfect exact matches}\\ \text{across the experimental groups}\end{pmatrix}. \quad \square$$

*Proof of Thm. 4.* Let $D_{ij} = \delta(X_i, X_j)$. The PSOD solves the optimization problem

$$\min_{W \in \mathcal{W}} \max_{\|f\|_{ip} \le 1} |B(W, f)| = \frac{2}{n}\min_{u \in \{-1,1\}^n : \sum_{i=1}^n u_i = 0} \max_{y \in \mathbb{R}^n : y_i - y_j \le D_{ij}} u^T y. \tag{E.1}$$

We will show that the set of optimal solutions $u$ to (E.1) is equal to the set of assignments of $+1, -1$ to the pairs in any minimal-weight pairwise match. Since the PSOD randomizes over these, this will show that it is equivalent to pairwise-matched allocation, which randomly splits pairs.

Consider any non-bipartite matching $\mu = \{\{i_1, j_1\}, \ldots, \{i_{n/2}, j_{n/2}\}\}$ and any $t \in \{-1, +1\}^{n/2}$. Let $u_{i_l} = t_l$, $u_{j_l} = -t_l$. Enforcing only a subset of the constraints on $y$, the cost of $u$ in (E.1) is bounded above as follows

$$\max_{y_i - y_j \le D_{ij}} u^T y = \max_{y_i - y_j \le D_{ij}} \sum_{l=1}^{n/2} t_l(y_{i_l} - y_{j_l}) \le \sum_{l=1}^{n/2} D_{i_l j_l},$$

which is the matching cost of $\mu$. Now let instead a feasible solution $u$ to (E.1) be given. Let $S = \{i : u_i = +1\} = \{i_1, \ldots, i_{n/2}\}$ and its complement $S^C = \{i : u_i = -1\} = \{i'_1, \ldots, i'_{n/2}\}$. By linear programming duality we have

$$\max_{y_i - y_j \le D_{ij}} u^T y = \min_{Fe - F^T e = u, F \ge 0} \sum_{i,j=1}^n D_{ij} F_{ij} \tag{E.2}$$

since the LHS is bounded $(\le D_{i_1 i'_1} + \cdots + D_{i_{n/2} i'_{n/2}})$ and feasible $(y_i = 0 \ \forall i)$. The RHS is an uncapacitated min-cost transportation problem with sources $S$ (with inputs 1) and sinks $S^C$ (with outputs 1). Consider any $j_s \in S$, $j_t \in S^C$ and any path $j_s, j_1, \ldots, j_p, j_t$. By the triangle inequality,

$$D_{j_s j_t} \le D_{j_s j_1} + D_{j_1 j_2} + \cdots + D_{j_p j_t}.$$

Therefore, it is always preferable to send flow along edges between $S$ and $S^C$ only. Thus, erasing all edges within $S$ or $S^C$, the problem is seen to be a bipartite matching problem. The min-weight bipartite matching is also a non-bipartite matching and by (E.2) its matching cost is the same as the cost of the given $u$ in the objective of (E.1). $\quad \square$

*Proof of Thm. 5.* The argument is similar to the above. This time the network flow problem has an additional node with zero external flow (neither sink nor source), uncapacitated edges into it from every other node with a unit cost of $\delta_0$, and uncapacitated edges out of it to every other node with a unit cost of $\delta_0$. $\quad \square$

*Proof of Thm. 6.* For the PSOD case we have,

$$
\begin{aligned}
M_P^2(W) &= \max_{k \ne k'} \max_{\|f\| \le 1} \left(\frac{1}{p}\sum_{i=1}^n (w_{ik} - w_{ik'})f(X_i)\right)^2 = \max_{k \ne k'} \left\|\left|\frac{1}{p}\sum_{i=1}^n (w_{ik} - w_{ik'})\mathcal{K}(X_i, \cdot)\right|\right\|^2\\
&= \frac{1}{p^2}\max_{k \ne k'} \left\langle \sum_{i=1}^n (w_{ik} - w_{ik'})\mathcal{K}(X_i, \cdot), \ \sum_{i=1}^n (w_{ik} - w_{ik'})\mathcal{K}(X_i, \cdot)\right\rangle\\
&= \frac{1}{p^2}\max_{k \ne k'} \sum_{i,j=1}^n (w_{ik} - w_{ik'})K_{ij}(w_{jk} - w_{jk'})
\end{aligned}
$$

Now, consider the maximum over $f$ in $M_M^2(Q)$. Let $f_0$ be a feasible solution. Write $f_0 = f + f^\perp$ with $f \in S = \text{span}\{\mathcal{K}(X_i, \cdot) : i = 1, \ldots, n\}$ and $f^\perp \in S^\perp$, its orthogonal complement. By orthogonality

$f^{\perp}(X_i) = \langle \mathcal{K}(X_i, \cdot), f^{\perp} \rangle = 0$ and $||f||^2 = ||f_0||^2 - ||f^{\perp}||^2 \leq 1$ so that $f$ achieves the same objective value as $f_0$ and remains feasible. Therefore we may restrict to $S$ and assume that $f = \sum_i \beta_i \mathcal{K}(X_i, \cdot)$ such that $\beta^T K \beta \leq 1$.

By positive semi-definiteness of $K$ and $Q$, we get

$$M_{\mathrm{M}}^2(Q) = \frac{2}{np} \sup_{\beta^T K \beta \leq 1} \sum_{i,j=1}^n Q_{ij} (K\beta)_i (K\beta)_j = \frac{2}{np} \sup_{\gamma^T \gamma \leq 1} \gamma^T \sqrt{K} Q \sqrt{K} \gamma = \frac{2}{np} \lambda_{\max} \left( \sqrt{K} Q \sqrt{K} \right). \quad \square$$

*Proof of Thm. 7.* By Asmp. 1(b), each $W_i$ by itself (but *not* the vector $W$) is statistically independent of $X, Y$ so that

$$\mathbb{E}[w_{ik} Y_{ik} \mid X, Y] = \mathbb{E}[w_{ik}] Y_{ik} = \tfrac{1}{m} Y_{ik}, \text{ and therefore}$$
$$\mathbb{E}[\hat{\tau}_{kk'} \mid X, Y] = \tfrac{1}{p} \sum_{i=1}^n \tfrac{1}{m} Y_{ik} - \tfrac{1}{p} \sum_{i=1}^n \tfrac{1}{m} Y_{ik'} = \mathrm{SATE}_{kk'}.$$

Note that we can rewrite $E_{kk'}$ as

$$E_{kk'} = \tfrac{1}{m} \sum_{l \neq k} \Xi_{kl} - \tfrac{1}{m} \sum_{l \neq k'} \Xi_{k'l} \quad \text{where} \quad \Xi_{kl} := \tfrac{1}{m} \sum_{i=1}^n (w_{ik} - w_{il}) \epsilon_{ik}.$$

Using the notation $A_{kk'} = \tfrac{1}{p} \sum_{i:W_i = k'} Y_{ik}$, $C_{kl} = B_{kl}(f_k) + \Xi_{kl}$, we have

$$\begin{aligned}
\hat{\tau}_{kk'} - \mathrm{SATE}_{kk'} &= A_{kk} - A_{k'k'} - \tfrac{1}{m} \sum_{l=1}^m A_{kl} + \tfrac{1}{m} \sum_{l=1}^m A_{k'l} \\
&= \tfrac{m-1}{m} A_{kk} - \tfrac{1}{m} A_{kk'} + \tfrac{1}{m} A_{k'k} - \tfrac{m-1}{m} A_{k'k'} \\
&\quad - \tfrac{1}{m} \sum_{l \neq k,k'} (A_{kk} - C_{kl}) + \tfrac{1}{m} \sum_{l \neq k,k'} (A_{k'k'} - C_{k'l}) = D_{kk'} + E_{kk'}.
\end{aligned}$$

Let $i, j$ be equal or unequal, $k, k', l, l'$ equal or unequal. Then,

$$\begin{aligned}
\mathrm{Cov}(w_{il} f_k(X_i), w_{jl'} \epsilon_{jk'}) &= \mathbb{E}[w_{il} w_{jl'} f_k(X_i) \mathbb{E}[\epsilon_{jk'} \mid X, Z]] \\
&\quad - \mathbb{E}[w_{il} f_k(X_i)] \mathbb{E}[w_{jl'} \mathbb{E}[\epsilon_{jk'} \mid X, Z]] = 0 - 0 = 0, \\
\mathrm{Cov}((w_{ik} - w_{il}) f_k(X_i), f_{k'}(X_j)) &= \mathbb{E}[w_{ik} - w_{il}] \mathrm{Cov}(f_k(X_i), f_{k'}(X_j)) = 0, \\
\mathrm{Cov}((w_{ik} - w_{il}) \epsilon_{ik}, f_{k'}(X_j)) &= \mathbb{E}[w_{ik} - w_{il}] \mathrm{Cov}(\epsilon_{ik}, f_{k'}(X_j)) = 0,
\end{aligned}$$

where the latter two equalities are due to the independence of $W_i$ due to blinding treatments. This proves uncorrelateness. The rest follows from an application of the law of total variance and rearranging terms. $\square$

*Proof of Thm. 8.* Define

$$Z(f,g) = \mathbb{E}\left[ \left( \tfrac{1}{p} \sum_{i=1}^n (w_{i1} - w_{i2}) f(X_i) \right) \left( \tfrac{1}{p} \sum_{i=1}^n (w_{i1} - w_{i2}) g(X_i) \right) \right].$$

By construction, $Z(f,f) \leq ||f||^2 \mathbb{E}[M_{\mathrm{opt}}^2]$. By Asmp. 1(b),

$$\begin{aligned}
\mathrm{Var}(B_{kl}(f)) &= Z(f,f) \quad \text{for } l \neq k, \\
\mathrm{Cov}(B_{kl}(f), B_{kl'}(f)) &= \tfrac{1}{2} Z(f,f) \quad \text{for } k, l, l' \text{ distinct,} \\
\mathrm{Cov}(B_{kl}(f), B_{k'l'}(g)) &= \begin{cases}
\tfrac{1}{2} Z(f,g) & \text{for } l = l' \notin \{k, k'\}, \\
-\tfrac{1}{2} Z(f,g) & \text{for } l = k', \, l' \neq k, \\
-\tfrac{1}{2} Z(f,g) & \text{for } l \neq k', \, l' = k, \\
-Z(f,g) & \text{for } l = k', \, l' = k, \\
0 & \text{for } k, k', l, l' \text{ distinct.}
\end{cases}
\end{aligned}$$

It follows that

$$
\begin{aligned}
\mathrm{Var}\,(D_{kk'}) &= \tfrac{1}{m^2}\left(\tfrac{m^2-m}{2}Z(f_k,f_k)+\tfrac{m^2-m}{m^2}Z(f_{k'},f_{k'})\right)+\tfrac{m+2}{m^2}Z(f_k,f_{k'})\\
&= \tfrac{1}{m^2}\left(\tfrac{m^2}{2}-m-1\right)\left(Z(f_k,f_k)+Z(f_{k'},f_{k'})\right)+\tfrac{1}{m^2}\left(\tfrac{m+2}{2}\right)Z(f_k+f_{k'},f_k+f_{k'})\\
&\leq \tfrac{1}{m^2}\left(\tfrac{m^2}{2}-m-1\right)\left(\mathbb{E}\left[M_{\mathrm{opt}}^2\right]||f_k||^2+\mathbb{E}\left[M_{\mathrm{opt}}^2\right]||f_{k'}||^2\right)+\tfrac{1}{m^2}\left(\tfrac{m+2}{2}\right)\mathbb{E}\left[M_{\mathrm{opt}}^2\right]||f_k+f_{k'}||^2\\
&\leq \tfrac{(||f_k||+||f_{k'}||)^2}{2}\left(1-\tfrac{1}{m}\right)\mathbb{E}\left[M_{\mathrm{opt}}^2\right]
\end{aligned}
$$

since $||f+g||^2\leq(||f||+||g||)^2$ and $\left(||f||^2+||g||^2\right)\leq(||f||+||g||)^2$. □

*Proof of Thm. 9.* This is immediate from Thm. 7 after observing that if $f_k$ is constant then $B_{kl}(W,f_k)=0$ for any $W$ and hence $D_{kk'}=0$. □

*Proof of Thm. 10.* Fix $f$ and $g$. Using $(\sum_{i=1}^b z_i)^2\leq b\sum_{i=1}^b z_i^2$,

$$
\begin{aligned}
Z(f,f) &= \mathbb{E}\Big(\tfrac{1}{p}\textstyle\sum_{i=1}^n (w_{i1}-w_{i2})(f-g)(X_i)+\tfrac{1}{p}\sum_{i=1}^n (w_{i1}-w_{i2})g(X_i)\Big)^2\\
&\leq 2\mathbb{E}\Big(\tfrac{1}{p}\textstyle\sum_{i=1}^n (w_{i1}-w_{i2})(f-g)(X_i)\Big)^2+2\mathbb{E}\Big(\tfrac{1}{p}\sum_{i=1}^n (w_{i1}-w_{i2})g(X_i)\Big)^2\\
&\leq \tfrac{2}{p^2}\times p\times p\times\tfrac{2}{m}\times\mathbb{E}((f-g)(X_1))^2+2Z(g,g)=\tfrac{4}{m}||f-g||_2^2+2Z(g,g)
\end{aligned}
$$

The rest is as in the proof of Thm. 8, choosing $g\in\mathcal{F}$. □

*Proof of Thm. 11.* First we prove that $\mathcal{F}$ is dense in $L^2$ with the measure of $X_1$. Simple functions $I_A(x)=\mathbb{I}[x\in A]$ for $A$ measurable are dense in $L^2$ (this is how one integrates). So it suffices to show that we can approximate $I_A(x)$. Fix any $\eta>0$. Let $U\supset A$ open and $E\subset A$ compact be such that $\mathbb{P}(X_1\in U)-\mathbb{P}(X_1\in E)\leq\eta/2$. By Urysohn's lemma (Royden, 1988), there exists a continuous function $f$ with support $K\subset U$ compact, $0\leq f\leq 1$, and $f(x)=1\,\forall x\in E$. Therefore, $\sqrt{\mathbb{E}[(I_A(X_1)-f)^2]}=\sqrt{\mathbb{E}[(I_A(X_1)-f)^2\mathbb{I}[X_1\in U\backslash E]]}\leq\mathbb{P}(X_1\in U\backslash E)\leq\eta/2$. Since $f$ is continuous with compact support, universality gives $\exists g\in\mathcal{F}$ such that $\sup_{x\in\mathcal{X}}|f(x)-g(x)|<\eta/2$. Since $\mathbb{E}[(f-g)^2]\leq\sup_{x\in\mathcal{X}}|f(x)-g(x)|^2$, we have $\sqrt{\mathbb{E}[(I_A(X_1)-g)^2]}\leq\eta$.

Let $\eta'=\sqrt{\eta/(8(m-1))}$. By the above, for each $\ell=k,k'$, $\exists g_\ell\in\mathcal{F}$ such that $||f_\ell-g_\ell||_2\leq\eta'$. The rest is immediate from Thm. 10. □

*Proof of Thm. 12.* Fix $W_i'=(i\bmod p)+1$. Define $\xi_i^{(k,k')}:f\mapsto\left(f\left(X_{m(i-1)+k}\right)-f\left(X_{m(i-1)+k'}\right)\right)$. Then, since $\xi_i^{(k,k')}$ is in the continuous dual space $\mathcal{F}^*$, we can write $M_{\mathrm{P}}^2(W')=\max_{k\neq k'}T_n^{(k,k')}$ where

$$
T_n^{(k,k')}=\sup_{||f||\leq 1}\left(\tfrac{1}{p}\textstyle\sum_{i=1}^p\xi_i^{(k,k')}(f)\right)^2=\left|\left|\tfrac{1}{p}\sum_{i=1}^p\xi_i^{(k,k')}\right|\right|_{\mathcal{F}^*}^2.
$$

Note $\xi_i^{(k,k')}$ are independent and identically distributed with expectation 0 (i.e., Bochner integral). $B$-convexity of $\mathcal{F}$ implies it has Rademacher type $>1$, which implies $\mathcal{F}^*$ has Rademacher type $>1$, which implies $\mathcal{F}^*$ is $B$-convex (Pisier, 2011). By $B$-convexity and the main result of Beck (1962) (or by Chen and Zhu, 2011 for the Hilbert case),

$$
T_n^{(k,k')}\to 0\text{ almost surely as }n\to\infty.
$$

As there are only finitely many $k,k'$, we have $M_{\mathrm{P}}^2(W')\to 0$ almost surely. By construction, $M_{\mathrm{M\text{-}opt}}^2\leq M_{\mathrm{P\text{-}opt}}^2\leq M_{\mathrm{P}}^2(W')$. Hence, the distance between $\hat{\tau}_{kk'}$ and $\mathrm{SATE}_{kk'}+E_{kk'}$ is $|D_{kk'}|\leq$ $\left(1 - \frac{1}{m}\right)\left(||f_k|| + ||f_{k'}||\right)\sqrt{M_{\text{opt}}^2} \to 0$ almost surely. Therefore, as $\text{SATE}_{kk'} + E_{kk'}$ is strongly consistent, so is $\hat{\tau}_{kk'}$.    □

*Proof of Thm. 13.* Let $\eta > 0$, $\delta > 0$ be given. By denseness proved in Thm. 11, there are $g_k, g_{k'} \in \mathcal{F}$ such that $||f_k - g_k||_2^2 \leq \delta\eta^2(1 - 1/m)^{-1}/32$. Let $\tilde{D}_{kk'} = \frac{1}{m}\sum_{l \neq k} B_{kl}(W, g_k) - \frac{1}{m}\sum_{l \neq k'} B_{k'l}(W, g_{k'})$. By Thm. 12, $\tilde{D}_{kk'} \to 0$ a.s., which implies $\tilde{D}_{kk'} \to 0$ in probability. Let $n_0$ be sufficiently large so that $\mathbb{P}\left(|\tilde{D}_{kk'}| > \eta/2\right) \leq \delta/2$ for all $n \geq n_0$. Since $\mathbb{E}[B_{kl}^2(W, f_k - g_k)] \leq 2||f_k - g_k||_2^2$, we have $\mathbb{E}(\tilde{D}_{kk'} - D_{kk'})^2 \leq 4\left(1 - \frac{1}{m}\right)||f_k - g_k||_2^2 \leq \delta\eta^2/8$. Therefore, for all $n \geq n_0$,

$$\mathbb{P}\left(|D_{kk'}| > \eta\right) \leq \mathbb{P}\left(|D_{kk'} - \tilde{D}_{kk'}| > \eta/2\right) + \mathbb{P}\left(|\tilde{D}_{kk'}| > \eta/2\right) \leq \frac{4}{\eta^2}\mathbb{E}(\tilde{D}_{kk'} - D_{kk'})^2 + \delta/2 \leq \delta.    □$$